\documentclass[12pt]{amsart}
\usepackage[dvips]{graphicx,color}
\usepackage[dvips]{epsfig}
\usepackage{amsmath,amssymb,latexsym,amscd,xypic}
\DeclareRobustCommand{\SkipTocEntry}[4]{}
\makeatletter
\newcommand\@dotsep{4.5}
\def\@tocline#1#2#3#4#5#6#7{\relax
  \ifnum #1>\c@tocdepth % then omit
  \else
    \par \addpenalty\@secpenalty\addvspace{#2}%
    \begingroup \hyphenpenalty\@M
    \@ifempty{#4}{%
      \@tempdima\csname r@tocindent\number#1\endcsname\relax
    }{%
      \@tempdima#4\relax
    }%
    \parindent\z@ \leftskip#3\relax \advance\leftskip\@tempdima\relax
    \rightskip\@pnumwidth plus1em \parfillskip-\@pnumwidth
    #5\leavevmode\hskip-\@tempdima #6\relax
    \leaders\hbox{$\m@th
      \mkern \@dotsep mu\hbox{.}\mkern \@dotsep mu$}\hfill
    \hbox to\@pnumwidth{\@tocpagenum{#7}}\par
    \nobreak
    \endgroup
  \fi}
\makeatother 

\DeclareFontFamily{OT1}{rsfs}{}
\DeclareFontShape{OT1}{rsfs}{n}{it}{<-> rsfs10}{}
\DeclareMathAlphabet{\curly}{OT1}{rsfs}{n}{it}

\newcommand\C{\mathbb C}
\newcommand\CC{\mathcal C}
\newcommand\FF{\mathbb F}

\newcommand\II{\mathbb I}
\newcommand\I{\curly I}

\newcommand\M{\mathcal M}

\renewcommand\O{\mathcal O}
\newcommand\PP{\mathbb P}

\newcommand\Q{\mathbb Q}

\newcommand\Z{\mathbb Z}

\newcommand{\rt}[1]{\stackrel{#1\,}{\rightarrow}}
\newcommand{\Rt}[1]{\stackrel{#1\,}{\longrightarrow}}
\newcommand\To{\longrightarrow}
\newcommand\into{\hookrightarrow}
\newcommand\Into{\ar@{^{ (}->}[r]}

\renewcommand\_{^{}_}
\newcommand\take{\backslash}
\newcommand\bull{{\scriptscriptstyle\bullet}}
\newcommand\udot{^\bull}

\newcommand\tr{\operatorname{tr}}

\newcommand\id{\operatorname{id}}

\newcommand\Hom{\operatorname{Hom}}
\renewcommand\hom{\curly H\!om}
\newcommand\Ext{\operatorname{Ext}}
\newcommand\ext{\curly E_{}xt}

\newcommand\Hilb{\operatorname{Hilb}}

\newcommand\Chow{\operatorname{Chow}}
\newcommand\Jac{\operatorname{Pic}}

\newcommand\beq[1]{\begin{equation}\label{#1}}
\newcommand\eeq{\end{equation}}
\newcommand\beqa{\begin{eqnarray*}}
\newcommand\eeqa{\end{eqnarray*}}

\makeatletter \@addtoreset{equation}{section} \makeatother

\newtheorem{defn}[equation]{Definition}
\newtheorem{thm}{Theorem}
\newtheorem{lem}[equation]{Lemma}

\newtheorem{prop}[equation]{Proposition}

\title[Stable pairs and BPS invariants]
{\textbf{Stable pairs and BPS invariants}}
\author{R. Pandharipande and R. P. Thomas}
\begin{document}

\begin{abstract} \noindent
We define the BPS invariants of Gopakumar-Vafa in the case of irreducible
curve classes on Calabi-Yau 3-folds. The main tools are
the theory of stable pairs in the derived category and Behrend's
constructible function approach to the virtual class. 
For irreducible curve classes, we prove
the stable
pairs generating function
satisfies the strong BPS rationality conjectures. 

We define the contribution
of each curve $C$ to the BPS invariants and show
the contributions lie 
between the geometric genus and arithmetic genus of $C$.
Complete formulae are derived for nonsingular
and nodal curves.

A discussion of primitive classes on $K3$ surfaces from
the point of view of stable pairs is given in the Appendix
via calculations of Kawai-Yoshioka.
A proof of the Yau-Zaslow formula for
rational curve counts is obtained. A connection
is made to the Katz-Klemm-Vafa formula for BPS counts
in all genera.
\end{abstract}

\maketitle \vspace{-5mm}
\thispagestyle{empty}
\tableofcontents

%%%%%%%%%%%%%%%%%%%%%%%%%%%%%%%%%%%%%%%%%%%%%%%%%%%%%%%%%%%%%%%%%%%%%%%%%%%

\setcounter{section}{-1}
\section{Introduction}
Let $X$ be a nonsingular, projective, Calabi-Yau 3-fold{\footnote{The Calabi-Yau condition
for us is $K_X\cong\O_X$. No restriction of the
fundamental group is necessary.}}.
Invariants counting curves in
$X$ via stable pairs have been defined in \cite{PT1}.
A {\em pair} $(F,s)$ consists of a sheaf $F$ on $X$  supported
in dimension 1 together with a section $s\in H^0(X,F)$.  A pair $(F,s)$
 is {\em stable} if 
\begin{enumerate}
\item[(i)]
the sheaf $F$ is {pure}, 
\item[(ii)] the section $\O_X \stackrel{s}{\rightarrow} F$ has 0-dimensional
cokernel.
\end{enumerate}
By purity (i), every nonzero subsheaf of $F$ has support of dimension 1.
As a consequence,
 the scheme theoretic support $C\subset X$ of  $F$ is a 
Cohen-Macaulay curve.
The support of the cokernel (ii) is a finite length subscheme $Z\subset C$.
If the support $C$ is nonsingular, then the stable pair $(F,s)$ is 
uniquely determined by $Z\subset C$.
However, for general $C$, the subscheme $Z$ does not determine $F$ and $s$.

Discrete invariants of a stable pair include the holomorphic Euler characteristic
$\chi(F)\in \mathbb{Z}$ and the class $[F]\in H_2(X,\mathbb{Z})$.
The moduli space $P_n(X,\beta)$ parameterizes stable pairs satisfying
$$\chi(F)=n, \ \ [F]=\beta.$$
Pair stability arises naturally as GIT stability in  the study of appropriate 
quotients
 \cite{LPPairs1,PT1}. In fact, GIT stability is found there 
to be equivalent to
semi-stability.  
The moduli 
space $P_n(X,\beta)$ is a therefore a projective scheme.

To define
invariants, we use a virtual cycle. The usual
deformation theory of pairs is not perfect in the sense of \cite{BFNormalCone}, but
the fixed-determinant deformation theory of the associated \emph{complex}
\beq{idot}
I\udot=\{\O_X\Rt{s}F\}\ \in\,D^b(X)
\eeq
is shown in \cite{PT1, HT} to define a perfect obstruction 
theory for $P_n(X,\beta)$ of virtual dimension zero. A virtual cycle 
is then obtained by \cite{BFNormalCone,LiTianVirtual}. The resulting
invariants
$$
P_{n,\beta}=\int_{[P_n(X,\beta)]^{vir}}1
$$
are conjecturally equal to the reduced DT invariants of \cite{MNOP1}.
Let
$$Z_\beta(q) = \sum_{n\in \mathbb{Z}} P_{n,\beta} \ q^n$$
be the generating series.
Calculations in the toric Calabi-Yau case can be found in \cite{PT2}.

Since $X$ is {Calabi-Yau},
the above deformation
theory of complexes is 
\emph{self-dual} in the sense of \cite{BehrendDT}. Heuristically,
$P_n(X,\beta)$ may be viewed  locally as the critical
locus of a function. 
The  virtual dimension of $P_n(X,\beta)$ is 0 and, on the nonsingular
locus of the moduli space, the obstruction sheaf is the cotangent bundle.
Therefore if $P_n(X,\beta)$ is everywhere nonsingular then
$$
P_{n,\beta}=(-1)^{\dim P_n(X,\beta)}\ e(P_n(X,\beta)),
$$
where $e$ denotes the topological Euler characteristic.

If singularities are present, $P_{n,\beta}$ certainly differs from 
the (signed) Euler characteristic.
By Behrend's results  \cite{BehrendDT}, there exists 
an \emph{integer-valued
constructible function} $\chi^B$ over any scheme with the property that if the scheme
is proper and admits a self-dual obstruction theory then the length of its virtual cycle
equals its $\chi^B$-weighted Euler characteristic. Therefore
$$P_{n,\beta}=e(P_n(X,\beta),\chi^B):=\sum_{n\in \Z} n\,e\big((\chi^{B})^{-1}(n)\big).
$$
At nonsingular points,
 $$\chi^B=(-1)^{\dim P_n(X,\beta)},$$
 but at singularities  $\chi^B$ is
a more complicated function.
The weighted Euler characteristic of $P_{n}(X,\beta)$
is a deformation invariant.

% 
% An {\em irreducible} class
% $\beta\in H_2(X,\beta)$ is a primitive and
% indecomposable class.
% All Cohen-Macaulay curves $C\subset H_2(X,\mathbb{Z})$ of
% irreducible class $\beta$
% are necessarily reduced and irreducible.
% Many of the ideas we present here apply
% to general classes $\beta\in H_2(X,\beta)$, but in the irreducible
% case there are fewer technical obstructions and our proofs
% are complete. Interesting irreducible classes can be found, for
% example, in 1-parameter families of $K3$ surfaces. 

Behrend's theory applied to $P_n(X,\beta)$ allows
us to use topological Euler characteristics and cut-and-paste techniques.
We require new technical results comparing the
value of Behrend's function at a pair $(F,s)$ to the value at the
sheaf $F$, see Theorem \ref{Kai} of Section \ref{one}. The arguments
turn out to be remarkably
simple when $\beta$ is an irreducible{\footnote{The
definitions of irreducible, primitive, and reduced curve
classes
are given in Section 1.}} class. 
We prove the following result in
Section \ref{BPSSec}. 

%The extension to reduced classes comes in Appendix
%\ref{reduced}.

\begin{thm} For $\beta$ irreducible, 
$Z_\beta(q)$ \label{aaa}
is the Laurent series expansion of a rational function in $q$.
\end{thm}

Serre duality relates a line bundle $L$ on a nonsingular curve
$C\subset X$ to $L^{-1} \otimes K_C$.
Since
$$\chi(L) = - \chi(L^{-1} \otimes K_C),$$
Serre duality relates the geometry of $P_n(X,\beta)$ to
$P_{-n}(X,\beta)$. The compatibility of
Serre duality 
with $\chi^B$ proven in  Section \ref{BPSSec} yields
a more subtle result.

\begin{thm} For $\beta$ irreducible,
the rational function $Z_\beta(q)$ \label{bbb}
is invariant under the transformation $q\leftrightarrow q^{-1}$.
\end{thm}

In fact, we prove $Z_\beta(q)$ satisfies the full BPS rationality conjectured
in \cite{PT1}.

\begin{thm} \label{ccc}
For $\beta$ irreducible,
\beq{bpsform}
Z_{\beta}(q)=\sum_{r=0}^gn_{r,\beta}\,q^{1-r}(1+q)^{2r-2}
\eeq
where the $n_{r,\beta}$ are integers and 
$g$ is the maximal arithmetic  
genus $g$ of a curve in class $\beta$.
\end{thm}

We obtain a deformation invariant definition
of the BPS counts $n_{r,\beta}$ of Gopakumar-Vafa \cite{GV1, GV2}
for irreducible classes $\beta$. 
In Section \ref{LocalSec}, we give a \emph{local} definition of these BPS
invariants for irreducible curve classes. We define
constructible functions over the space of curves in $X$ with respect to which
the weighted Euler characteristics yields the BPS numbers.
We prove the functions  
are nonzero on $C\subset X$ only in genus $g$ between the geometric
and arithmetic genera of $C$. 
Complete evaluations of the functions 
are obtained for  nonsingular and nodal curves.

In Appendix A, we sketch the extension of Theorems
1 and 2 to reduced curve classes which are not necessarily irreducible. We
also explain what is needed to show the
vanishing of 
BPS counts in negative genus in the reduced case.

Interesting examples of irreducible and reduced classes
occur on surfaces.
If $C$ is Gorenstein, the stable pairs with support $C$ are
proven in Appendix \ref{B} to correspond  
bijectively to finite length subschemes $Z\subset C$.
The moduli spaces of stable pairs {\em on a surface} are
then shown to be isomorphic to relative Hilbert schemes.

In Appendix C, the beautiful theory of primitive classes
on $K3$ surfaces is considered. By results of Kawai-Yoshioka \cite{ky}, 
the Katz-Klemm-Vafa \cite{KKVSpinning} formula for BPS state counts 
is obtained for the
theory of stable pairs. The corresponding calculations in
Gromov-Witten theory have not yet been completed.{\footnote{A discussion
of the Gromov-Witten side can be found in \cite{dm}.}

Let
$r_{0,g}$ be the number of rational curves of fixed primitive
class
with self-intersection $2g-2$ on a $K3$ surface.
Using the genus 0 BPS counts
together with the local BPS
theory of Section \ref{LocalSec},  a new proof of the
Yau-Zaslow formula,
$$\sum_{g\geq 0} r_{0,g} q^g = \prod_{n\geq 0} (1-q^n)^{-24},$$
is obtained.

The Yau-Zaslow formula was  proven in the primitive{\footnote{A
proof of the Yau-Zaslow formula for all curve classes on $K3$
surfaces has been recently obtained in \cite{kmps}.}} case
by Bryan-Leung \cite{bl} via Gromov-Witten theory. Our proof is very
close in spirit to the original sheaf-theoretic motivations for the formula
\cite{yz}. In particular, our argument via stable pairs and BPS
counts is parallel to Beauville's proof using compactified Jacobians
and Euler characteristics \cite{bea}.

\addtocontents{toc}{\SkipTocEntry}
\subsection*{Acknowledgements}  
We thank J. Bryan and
D. Maulik for many conversations
related to stable pairs and
A. Marian for pointing out the basic connection to
the work of Kawai-Yoshioka. 
We are grateful to S. Kleiman for advice on Jacobians of singular curves,
and Hua-Liang Chang for a careful
reading of the manuscript.

R.P. was partially supported by NSF grant DMS-0500187 and a Packard foundation
fellowship. R.T. was partially supported
by a Royal Society University Research Fellowship.
He also thanks the Leverhulme Trust and Columbia University for
a visit to New York in the spring of 2007 when the project was started. 
Many of the results presented here were found during a visit
to Lisbon in the summer of 2007.

\section{$\chi^B$-functions} \label{one}
Let $X$ be a nonsingular projective variety over $\C$. A nonzero
class $\beta\in H_2(X,\Z)$
is
{\em effective} if $\beta$ is represented by an algebraic curve.

\begin{defn} \label{superprimitive}
An effective class $\beta\in H_2(X,\Z)$ is 

\vspace{6pt}

\noindent $\bullet$ \emph{irreducible} if there
is no decomposition $\beta=\beta_1+\beta_2$ into
nonzero effective classes $\beta_i$,

\vspace{6pt}

\noindent $\bullet$  \emph{primitive} if $\beta$ is not a positive integer multiple of an effective
class,

\vspace{6pt}

\noindent $\bullet$  \emph{reduced} if in every
 decomposition $\beta=\sum_i\beta_i$ into effective
classes, all of the $\beta_i$ are primitive.

\end{defn}

For example, classes $\beta$ of minimal degree $\int_\beta c_1(L)$ measured
against any ample class $L$ are irreducible. 
Any primitive class $\beta$ on a $K3$ surface $S$ is irreducible on a generic
deformation of $S$ for which $\beta$ is of type $(1,1)$.

Let $X$ be a Calabi-Yau 3-fold.        
If $(F,s)$ is a
 a stable pair of irreducible class $\beta$,
$F$ is a {\em stable sheaf} 
since all quotient sheaves  have 0-dimensional
support.{\footnote{Usually, sheaf stability depends
upon the polarization $L$. However, for sheaves $F$ of rank 1 on an
 irreducible curve in $X$, stability is equivalent to purity. 
And, no strictly semi-stable sheaves exist.}} 
There is therefore a map
\beq{themap}
P_n(X,\beta)\Rt{\phi_n}\M_n(X,\beta)
\eeq
from the moduli space of stable pairs to the moduli space of stable pure
sheaves of Hilbert polynomial 
$$\chi(F(k))=k\int_\beta c_1(L)+n.$$
Moreover, the fibre of \eqref{themap} over a point $\{F\}$ is $\PP(H^0(F))$.
By the irreducibility of $\beta$, the cokernel of any nonzero $s$ section
is $0$-dimensional and $(F,s)$ is a stable pair.

Since $X$ is Calabi-Yau, both $P_n(X,\beta)$
\cite{PT1} and
$\M_n(X,\beta)$ \cite{ThCasson} 
  have self-dual obstruction theories. We can therefore
apply the results of \cite{BehrendDT}.

\begin{lem} The obstruction theory of $P_n(X,\beta)$ obtained
from fixed determinant deformations in the derived category \cite{PT1} is
self-dual in the sense of Behrend \cite{BehrendDT}.
\end{lem}

\begin{proof}
The obstruction theory of \cite{PT1} can be described as follows. Let
$$
\pi\colon X\times P_n(X,\beta)\to P_n(X,\beta)
$$
be the projection. There is a universal stable pair \cite{PT1},
$$
\O_{X\times P_n(X,\beta)}\to\FF,
$$
over $X\times P_n(X,\beta)$. Let $\II\udot$ be the associated
complex (with $\O_{X\times P_n(X,\beta)}$ in degree 0). Consider the complex
\begin{equation}\label{zz2}
R\pi_*R\hom(\II\udot,\II\udot\otimes\omega_\pi)_0[2]
\end{equation}
of trace-free Exts, where $\omega_\pi$ denotes the relative canonical bundle.
In \cite{PT1}, the complex \eqref{zz2}
is shown to be quasi-isomorphic to a 2-term complex
of locally free sheaves $\{E_1\to E_0\}$ over $P_n(X,\beta)$, with a canonical
morphism
$$
R\pi_*R\hom(\II\udot,\II\udot\otimes\omega_\pi)_0[2]\to L\udot_{P_n(X,\beta)}
$$
to the cotangent complex of $P_n(X,\beta)$.
The morphism is
 obstruction
theory for $P_n(X,\beta)$: the induced maps on $h^0$ and $h^1$
 are  isomorphisms and surjections respectively.

For $X$ Calabi-Yau, $\omega_\pi$ is trivial. Therefore,
 by relative Serre duality
for $\pi$, we obtain a quasi-isomorphism
$$
R\pi_*R\hom(\II\udot,\II\udot)_0^\vee\,\simeq\,
R\pi_*R\hom(\II\udot,\II\udot)_0[3].
$$
Thus
$$
\{E_0^\vee\to E_1^\vee\}[1]\ \simeq\ \{E_1\to E_0\},
$$
which is the definition of self-duality in  \cite{BehrendDT}.
\end{proof}

For any scheme $M$, Kai Behrend  \cite{BehrendDT} defines a canonical constructible
function 
$$\chi^B: M\rightarrow \Z,$$ depending only on the local scheme 
structure{\footnote{In fact,
$\chi^B$ depends only on the
local scheme structure  in analytic topology by Proposition 4.22 of \cite{BehrendDT}.}}.
If $M$ is compact and equipped with a self-dual obstruction theory, then
$$\int_{[M]^{vir}}1 = e(M,\chi^B)$$
where the right side is the weighted Euler characteristic
$$e(M,\chi^B)= \sum_{n\in \Z} n\  e(\,(\chi^{B})^{-1}(n)\,)$$
and $e$ is the usual topological Euler characteristic.

If $M$ is nonsingular, then $\chi^B$ is 
the constant function $(-1)^{\dim M}$ and 
$$\int_{[M]^{vir}} 1 = (-1)^{\dim M}e(M).$$
 More generally, by  Proposition 1.5(i) of \cite{BehrendDT},
if $f\colon M\to N$ is a \emph{smooth} map of relative dimension $r$, then
\beq{smooth}
\chi^B_M=(-1)^rf^*\chi^B_N.
\eeq
If $e(F(f))$ is the Euler characteristic of the fibre of $f$, then
$$e(M,\chi^B_M)= (-1)^re(N,\chi^B_N)\cdot e(F(f)).$$

On $P_n(X,\beta)$ and $\M_n(X,\beta)$,
 we obtain functions $\chi\_P$
and $\chi\_\M$.{\footnote{For convenience, we will often drop the
superscripted $B$ in the notation.}} 
The invariants
\begin{eqnarray*}
P_{n,\beta} =\int_{[P_n(X,\beta)]^{vir}} 1& =& e(P_n(X,\beta), \chi\_P),\\  
N_{n,\beta} =\int_{[\M_n(X,\beta)]^{vir}}1& =& e(\M_n(X,\beta),\chi\_\M)
\end{eqnarray*}
are the weighted Euler characteristics.

The following property holds even though
 the map \eqref{themap} may be neither smooth nor surjective.
The result underpins the whole paper.

\begin{thm} \label{Kai}
$\chi\_P=(-1)^{n-1}\phi_n^*\chi\_\M$.
\end{thm}

\begin{proof}
We work locally around one point $(F,s)$ of $P_n(X,\beta)$. 
By the irreducibility
of $\beta$, the Cohen-Macaulay support $C$ of $F$ is reduced and irreducible.
Hence, $C$ is generically nonsingular. There exists 
a local smooth divisor $D_C$ which intersects $C$ (and all nearby $C$ in
the same homology class) transversally in a single point. We may also assume
$D_C\cap C$ to be disjoint from the zeros of $s$.

Let $k\ge0$. Tensoring with $\O(kD_C)$ and the canonical section $s_{D_C}^k$
yields a map of analytic open sets:
\beq{miracle} \xymatrix{
P_n(X,\beta)\supset\hspace{-1cm} & V_n \rto\dto^{\phi_n} & V_{n+k} \dto^{\phi_{n+k}}
& \hspace{-1cm} \subset P_{n+k}(X,\beta) \\
\M_n(X,\beta)\supset\hspace{-8mm} & U_n \rto^(.4)\simeq & U_{n+k} & \hspace{-8mm}
\subset\M_{n+k}(X,\beta)}\ .
\eeq
Here, $U_n$ is a sufficiently small analytic neighborhood of $F$.
Since
$\chi^B$-functions depend only on the local scheme structure
\cite{BehrendDT}, the bottom isomorphism makes the $\chi^B$-functions 
of the two
sheaf moduli spaces $\M$ locally the same. We call them $\chi\_\M$.
The open sets
$$V_n=\phi_n^{-1}(U_n), \ \ V_{n+k}=\phi_{n+k}^{-1}(U_{n+k})$$
contain $(F,s)$ and $(F(kD_C),s.s_{D_C}^k)$ respectively.
For $k$ sufficiently large, $\phi_{n+k}$ is a smooth $\PP^{n+k-1}$-bundle.
 
By making $V_{n+k}$ smaller if necessary, the map $V_n\to V_{n+k}$ admits
a left inverse
$$
V_n\stackrel{\,\psi}{\leftarrow}V_{n+k}
$$
given by forgetting about the $k$ points close to $D_C$. The map $\psi$
is locally 
smooth  with fibre the $k$th symmetric product of an open subset of
a nonsingular curve.

We calculate the $\chi^B$-function of $V_{n+k}$
in two different ways round the commutative diagram \eqref{miracle}, using
\eqref{smooth} applied to the two smooth maps $\psi$ and $\phi_{n+k}$.
The two resulting expressions are
$$
(-1)^{n+k-1}\phi_{n+k}^*\chi\_\M\ =\ (-1)^k\psi^*\chi\_P.
$$
Pulling back to $V_n$ gives
$$
(-1)^{n+k-1}\phi_n^*\chi\_\M\ =\ (-1)^k\chi\_P.
$$
Multiplying by $(-1)^k$ gives the result. \end{proof}

\section{BPS rationality} \label{BPSSec}
\subsection{Results}
Let $X$ be a Calabi-Yau 3-fold and $\beta\in H_2(X,\Z)$
an irreducible class. Let $g$ be the maximal arithmetic genus of a curve
in the class $\beta$. Following the notation of Section \ref{one}, let
$$P_{n,\beta}=\int_{[P_n(X,\beta)]} 1, \qquad
N_{n,\beta} = \int_{[\M_n(X,\beta)]^{vir}} 1$$
denote the invariants of \cite{PT1,ThCasson}.

\begin{prop} 
For $\beta$ an irreducible class, $N_{1,\beta}=N_{n,\beta}$
for all $n$.
\end{prop}

\begin{proof}
Let $F$ be a stable sheaf determining a moduli point of
$\M_1(X,\beta)$. Let $C$ be the support of $F$. As
in the proof of Theorem \ref{Kai}, let $D_C$ be a transverse divisor
meeting $C$ in 1 point. Let
$U_1\subset \M_1(X,\beta)$ be the analytic open set 
of sheaves supported on curves with a single transverse intersection
with $D_C$. 
%The open set $U_1$ is stable under the following operation:
%if ${F'}\in U_1$ and $L'$ is a degree 0 line bundle on the support
%$C'$ of $F'$, then ${F'\otimes L'}\in U_1$.
Tensoring with multiples of $\O(D_C)$
makes $U_1$ isomorphic to a corresponding
open set $U_n$ in each $\M_n(X,\beta)$. 

If $\M_1(X,\beta)$ is
covered by finitely many open sets $U^i_1$ of the above form, the 
corresponding open
sets $U^i_n$ cover $\M_n(X,\beta)$.
By construction, the intersections
$$U^i_1 \cap U^j_1,\ \ U^i_1\cap U^j_1 \cap U^k_1,\ \ldots$$
are isomorphic to the corresponding intersections
$$U^i_n \cap U^j_n,\ \ U^i_n\cap U^j_n \cap U^k_n,\ \ldots \ .$$
 Calculating the weighted
Euler characteristics of the spaces $\M_1(X,\beta)$ and $\M_n(X,\beta)$
 as a sum of weighted Euler characteristics of the $U^i$
(minus the weighted Euler characteristics of their intersections,
plus the triple intersections and so on), we find 
$
N_{1,\beta}=N_{n,\beta}$.
\end{proof}

\begin{prop}\label{ggtc} Let $\beta\in H_2(X,\Z)$ be an irreducible class.
The invariants $P_{n,\beta}$  satisfy the following identities:
% $$
% \begin{array}{rclcl}
% P_{n,\beta} &=& (-1)^{n-1}n\,N_{1,\beta}, &\qquad& g\le n, \label{g0} \\
% P_{n,\beta}-P_{-n,\beta} &=& (-1)^{n-1}n\,N_{1,\beta}, && -g<n<g. \label{gg}
% \\
% P_{n,\beta} &=& 0, && n\le-g. \label{0}
% \end{array}
% $$
\begin{eqnarray}
P_{n,\beta}\ =\ (-1)^{n-1}n\,N_{1,\beta}, &\qquad& g\le n, \label{g0} \\
P_{n,\beta}-P_{-n,\beta}\ =\ (-1)^{n-1}n\,N_{1,\beta}, && -g<n<g. \label{gg}
\\
P_{n,\beta}\ =\ 0, \hspace{23mm} && n\le-g. \label{0}
\end{eqnarray}
\end{prop}

\begin{proof}
An element $(F,s)$ of $P_n(X,\beta)$ yields an exact sequence
$$
0\to\O_C\stackrel{s}{\to} F\to Q\to0,
$$
where $Q$ has 0-dimensional support. We obtain the inequality
$$
n=\chi(F)=\chi(\O_C)+\chi(Q)\ge 1-g+0>-g.
$$
Therefore $P_n(X,\beta)$ is empty for $n\le-g$, which implies \eqref{0}.
We verify \eqref{g0} and \eqref{gg}
simultaneously by proving \eqref{gg} \emph{for all} $n\in\Z$.

If $F$ is a line bundle on a nonsingular curve $C\subset X$, then 
Serre duality relates $F$ and $F^{-1}\otimes K_C$. More generally,
there is a map
\begin{eqnarray}
\M_n(X,\beta) &\to& \M_{-n}(X,\beta) \nonumber \\
F &\mapsto& \ext^2_X(F,K_X). \label{duall}
\end{eqnarray}
Since $F$ is pure, $F$ has homological dimension 2 \cite[Proposition
1.1.10]{HLShaves} so  $\ext^{\ge3}(F,K_X)=0$. Similarly $\ext^{\le1}(F,K_X)=0$
because $F$ is supported in codimension 2. Therefore
$$
\ext^2(F,K_X)\cong R\hom(F,K_X)[2],
$$
which has the same Chern classes as elements of $\M_{-n}(X,\beta)$.

Pick a 3-term locally free resolution of $F$,
$$
0\to F_2\to F_1\to F_0\to F\to0.
$$
Applying $\hom(\ \cdot\ ,K_X)$ gives a 3-term locally free resolution
$$
0\to F_0^*\otimes K_X\to F_1^*\otimes K_X\to F_2^*\otimes K_X\to\ext^2(F,K_X)\to0
$$
of $\ext^2(F,K_X)$.
Therefore by \cite[Proposition 1.1.10]{HLShaves} $\ext^2(F,K_X)$ is
a pure sheaf. By the irreducibility assumption, $\ext^2(F,K_X)$
is stable and indeed defines an element of $\M_{-n}(X,\beta)$.

The map \eqref{duall}
is an involution and hence yields an isomorphism
$$
\M_n(X,\beta)\cong\M_{-n}(X,\beta).
$$
We may therefore consider the projections $\phi_n$ and $\phi_{-n}$
\eqref{themap} to fibre $P_n(X,\beta)$ and $P_{-n}(X,\beta)$ over the \emph{same}
space $\M_{n,\beta}$. We have
$$
H^0(\ext^2(F,K_X))\cong\Ext^2(F,K_X)\cong H^1(F)^*
$$
by Serre duality on $X$.
 The fibres of $\phi_n$ and $\phi_{-n}$ over $F\in\M_n(X,\beta)$
are therefore
\beq{fibres}
\PP(H^0(F)) \quad\text{and}\quad \PP(H^1(F)^*)
\eeq
respectively.

We stratify $\M_n(X,\beta)$ by the dimension of $H^0(F)$,
$$\M_n(X,\beta)=\cup_r V_r,$$
where $V_r$ is
the locus of sheaves $F$ with $h^0(F)=r$. There
are induced stratifications $\phi_{\pm
n}^{-1}(V_r)$ of $P_n(X,\beta)$ and $P_{-n}(X,\beta)$.
By \cite{BehrendDT}, we may calculate the invariants $P_{\pm n,\beta}$ via
these stratification as
$$
P_{\pm n,\beta}=\sum_re\big(\phi^{-1}_{\pm n}(V_r),\chi\_P|\_{\phi^{-1}_{\pm
n}(V_r)}\big)=\sum_r(-1)^{n-1}e(\phi^{-1}_{\pm n}(V_r),\phi_{\pm n}^*\chi\_\M),
$$
with the last equality following from Theorem \ref{Kai}. 
The $\chi^B$-function is the \emph{same constant} on the
fibres of both fibrations $\phi_{\pm n}$.

By \eqref{fibres}, over $V_r$, $\phi_n$ is a $\PP^{r-1}$-bundle and $\phi_{-n}$
is a $\PP^{r-n-1}$-bundle. These fibres have Euler characteristics $r$ and
$r-n$ respectively. We find
$$
P_{n,\beta}=\sum_r(-1)^{n-1}r\,e(V_r,\chi\_\M)
$$
and
$$
P_{-n,\beta}=\sum_r(-1)^{n-1}(r-n)e(V_r,\chi\_\M).
$$
Subtracting gives
\begin{eqnarray*}
P_{n,\beta}-P_{-n,\beta} &=& (-1)^{n-1}n\sum_re(V_r,\chi\_\M) \\
&=& (-1)^{n-1}n\,e(\M_n(X,\beta),\chi\_\M)=(-1)^{n-1}n\,N_{n,\beta}.
\end{eqnarray*}
\end{proof}

By Proposition \ref{ggtc}, the 
generating series 
\beq{Lseries}
Z_{\beta}(q)=\sum_nP_{n,\beta}\ q^n
\eeq
is the Laurent expansion of rational function in $q$, completing
the proof of Theorem \ref{aaa}.
However, a stronger statement can be made. 
Any Laurent series such as \eqref{Lseries} can be written as
\beq{PT1BPS}
Z_{\beta}(q)=\sum_rn_{r,\beta}q^{1-r}(1+q)^{2r-2},
\eeq
where the sum is over all $r\in\Z$ and only finitely many terms
with $r\ge0$ are nonzero; see \cite{PT1}. 
Moreover the integrality of the coefficients $P_{n,\beta}$ of
\eqref{Lseries} is
equivalent to the integrality of the $n_{r,\beta}$.

The conditions (\ref{g0}--\ref{0}) easily
imply the vanishing of $n_{r,\beta}$ for
$r<0$ and $r>g$. Therefore, by Proposition \ref{ggtc}, $Z_{P,\beta}(q)$
can be written uniquely in the BPS form
\beq{BPS}
Z_{\beta}(q)=\sum_{r=0}^gn_{r,\beta}\,q^{1-r}(1+q)^{2r-2} 
\eeq
for \emph{integers} $n_{r,\beta}$ which vanish for $r<0$ and for $r$ greater
than the largest genus $g$ of a holomorphic curve in the class $\beta$.
Since \eqref{BPS} is invariant under 
$q\leftrightarrow q^{-1}$, Theorems \ref{bbb} and \ref{ccc} are
proven.

\subsection{Remarks}
From formula \eqref{BPS}, 
we find the genus $0$ BPS invariant $n_{0,\beta}$ equals
is the DT invariant of sheaves $N_{1,\beta}$ 
in agreement with the proposal of
S. Katz \cite{KatzBPS}.
In fact, Katz expects $$n_{0,\beta}=N_{1,\beta}$$
to hold in much greater generality.
\medskip

Our sheaf theoretic definition of BPS invariants \eqref{BPS} in the irreducible
case is
the first rigorous and manifestly deformation-invariant approach.
Other papers on
the subject \cite{HST, SSGV, TodaGV} have defined BPS invariants
following the original perspective
of \cite{GV1, GV2} using $(\mathfrak{sl}_2\times\mathfrak{sl}_2)$-actions
on sophisticated
cohomology theories, but have been unable
to incorporate the virtual class. These definitions are therefore unlikely
to be deformation
invariant. Our definition is rather simpler, and more in line with the viewpoint
of \cite{KKVSpinning}.\medskip

It should be possible to extend our results to the Fano case
for \emph{any}
class $\beta$. After imposing
the requisite number of incidence conditions to cut the virtual dimension
to 0, the Fano case behaves like the Calabi-Yau case for \emph{irreducible}
$\beta$, as all other invariants vanish.
However, at present,
the analogue of $\chi^B$ is missing in the Fano case. \medskip

%As an example, one can impose the incidence conditions that the curves all
%pass through the points $D_0\cap C_0$, where $D_0$ is a fixed anticanonical
%divisor and $C_0$ a fixed curve in the class $\beta$. The deformation theory
%of such curves can be shown to be self-dual, with, in the
%smooth case, tangents governed by
%$$
%H^0(\nu\_C(-D_0))=H^0(\nu\_C\otimes K_X)\cong H^0(\nu_C^*\otimes K_C),
%$$
%and obstructions by its Serre dual $H^1(\nu\_C)$. The above results
%can probably then be extended to this case. For other incidence conditions,
%however, one must look for a different technique.

\subsection{Wall-crossing} \label{abayer}
Arend Bayer \cite{BayerPoly} and Yukinoba Toda \cite{TodaStab} have made
the beautiful observation that \eqref{gg} should be seen
as a wall-crossing formula. In fact, the wall-crossing
is much simpler than the wall-crossing
 conjectured in \cite{PT1} to equate the invariants $P_{n,\beta}$
to the reduced DT invariants of \cite{MNOP1}. For any 
$$I\udot=\{\O_X\to F\}\in
P_n(X,\beta),$$ we have the obvious exact triangle
\beq{stab}
F[-1]\to I\udot\to\O_X.
\eeq
Taking the derived dual gives
\beq{unstab}
\O_X\to(I\udot)^\vee\to\ext^2(F,K_X)[-1],
\eeq
where $\ext^2(F,K_X)\in\M_{-n}(X,\beta)$ is the sheaf dual to $F$ under the
duality \eqref{duall}.

Start with a stability
condition for which the complexes $I\udot\in P_{n,\beta}$ and the sheaves
$F,\,\O_X$ are stable. In
particular, the phase of $F[-1]$ should be less than that of $\O_X$ due to
the exact triangle \eqref{stab}.
Now pass through a codimension 1 wall in the space of stability conditions
so that the
phase of $F[-1]$ crosses that of $\O_X$. The extensions \eqref{stab}
become unstable, while extensions in the opposite direction \eqref{unstab}
become stable. Therefore, on the other side of the wall,
the stable objects are the derived duals of the complexes
made out of stable pairs in $P_{-n}(X,\beta)$.

Ideally, wall-crossing should be studied with Bridgeland stability conditions.
However, at present, their existence is conjectural. If instead we use Bayer's polynomial
stability conditions or Toda's limit stability conditions, then the analysis
can be made precise. These
stability conditions have been constructed, and the stable objects are as
claimed above \cite{BayerPoly, TodaStab}.

Since the pieces $\O_X,\,F,\,\ext^2(F,K_X)$ occurring in the  complexes are
also stable in these stability conditions, Joyce's conjectural
wall-crossing formula \cite{Joy} takes a very simple form. 
We count
only complexes of trivial determinant throughout. The invariant counting
the stable objects $I\udot\in P_n(X,\beta)$ on one side of the wall should
differ from those $(I\udot)^\vee\in P_{-n}(X,\beta)$ on the other side by
\begin{equation}\label{jjj2}
(-1)^{\chi(\O_X,F[-1])}\chi(\O_X,F[-1])\cdot\#(\O_X)\cdot\#(F),
\end{equation}
where 
$$\chi(\O_X,F[-1])=\sum_i(-1)^i\dim\Ext^i(\O_X,F[-1])=-n,$$
 and $\#$ denotes
the virtual number of elements of the moduli space of stable objects of the
corresponding type. 
For us, \eqref{jjj2} predicts
\begin{equation}\label{ntre}
P_{n,\beta}-P_{-n,\beta}=(-1)^{-n}(-n)\cdot 1\cdot N_{n,\beta}
\end{equation}
in precise agreement with \eqref{gg}. Perhaps \eqref{ntre} is
the first nontrivial example
of a wall-crossing formula in the derived category that can be rigorously
proved. \medskip

Toda \cite{TodaStab} has gone further with  wall crossings for arbitrary
(rather than irreducible) stable pairs. Using the
work of Joyce \cite{Joy}, he proves analogues of Theorems \ref{aaa}
and \ref{bbb}
for the Euler characteristics of the moduli spaces of stable pairs.
Once Behrend's function $\chi^B$ and the identities of Kontsevich-Soibelman
about the value $\chi^B$ on extensions \cite{KS} have been incorporated into
Joyce's work, Theorems \ref{aaa}
and \ref{bbb} for all classes on Calabi-Yau 3-folds should
be obtained.

\section{Local definition of BPS invariants} \label{LocalSec}
\subsection{Fixed curve}
Let $X$ be a Calabi-Yau 3-fold.
Throughout this Section we fix a Cohen-Macaulay curve $C\subset X$ in the
irreducible class $\beta$ of arithmetic genus $g=1-\chi(\O_C)$. 
The curve $C$ is reduced and irreducible. Let
$$
P_n(C)\subset P_n(X,\beta) \quad\text{and}\quad \M_n(C)\subset\M_n(X,\beta)
$$
denote the loci of stable pairs and pure sheaves supported on
$C$. Define 
 localised invariants\footnote{These 
are not invariants of $C$ alone. The dependence
on the embedding $C\subset X$ comes through $\chi\_P$.} by
\beq{PnC}
P_{n,C}=e\left(P_n(X,C),\chi\_P|\_{P_n(C)}\right).
\eeq

In Proposition \ref{ggtc}, we computed the
weighted Euler characteristics
of the spaces $P_n(X,\beta)$
using the map \eqref{themap} to $\M_n(X,\beta)$.
We can instead restrict attention
to the loci $P_n(C)\subset P_n(X,\beta)$, the inverse images of $\M_n(C)\subset
\M_n(X,\beta)$.
 The same proof applies, since $\bigcup_n\M_n(C)$ is invariant
under tensoring by line bundles and under the duality \eqref{duall}. 
We therefore obtain the same identities for localised invariants:
$$
\begin{array}{rclcl}
P_{n,C} &=& (-1)^{n-1}n\,N_{1,C}, &\qquad& g\le n, \\
P_{n,C}-P_{-n,C} &=& (-1)^{n-1}n\,N_{1,C}, && -g<n<g. \\
P_{n,C} &=& 0, && n\le-g,
\end{array}
$$
where $N_{1,C}=e\left(\M_n(C),\chi\_\M|\_{\M_n(C)}\right)$. 
Thus, the generating
series
\beq{ZPC}
Z_{C}(q)=\sum_n P_{n,C}\ q^n
\eeq
can be written uniquely as
\beq{localBPS}
Z_{C}(q)=\sum_{r=0}^gn_{r,C}\,q^{1-r}(1+q)^{2r-2} 
\eeq
for integers $n_{r,C},\ r=0,\ldots,g$.
\medskip

\subsection{Chow}
Let $\Chow(X,\beta)$ denote the variety of 1-dimensional cycles in the
class $\beta$. Since $\beta$ is irreducible, the cycles have
no multiplicities.
In fact, $\Chow(X,\beta)$ 
parameterises Cohen-Macaulay curves in class $\beta$.

The spaces $P_n(X,\beta)$ and $\M_n(X,\beta)$ map to
$\Chow(X,\beta)$ with fibres $P_n(C)$ and $\M_n(C)$ respectively.
We may calculate weighted Euler characteristics of $P_n(X,\beta)$ and
$\M_n(X,\beta)$ as weighted Euler characteristics of $\Chow(X,\beta)$, with
weight function the weighted Euler
characteristics of the fibres. 
More precisely, the integers $P_{n,C}$ \eqref{PnC}
define constructible functions
\begin{eqnarray*}
\Chow(X,\beta) &\to& \Z, \\
C &\mapsto& P_{n,C\,},
\end{eqnarray*}
whose weighted Euler characteristics are the
integers $P_{n,\beta}$.
 Similarly, $Z_{C}(q)$ \eqref{ZPC}
defines a $\Z(\!(q)\!)$-valued constructible function on $\Chow(X,\beta)$
with weighted Euler characteristic $Z_{\beta}(q)$.

Therefore the $n_{r,C}$ \eqref{localBPS} define constructible functions
\begin{eqnarray*}
\tilde n_{r,\beta}\colon \Chow(X,\beta) &\to& \Z, \\
C &\mapsto& n_{r,C\,},
\end{eqnarray*}
such that the BPS invariants $n_{r,\beta}$ of \eqref{BPS} are the weighted
Euler characteristics
$$
n_{r,\beta}=e(\Chow(X,\beta),\tilde n_{r,\beta}).
$$
We call $\tilde n_{r,\beta}(C)=n_{r,C}$ the \emph{contribution of $C\subset
X$ to $n_{r,\beta}$}.
\bigskip

Since these definitions hide behind a lot of formulae, their
exact meaning is
rather opaque. We would like to be able to compute the 
contributions $n_{r,C}$
directly, without computing all of the stable pairs invariants $P_{n,C}$.
In fact, the invariants $n_{r,\beta}$ are the more fundamental invariants,
from which all others (GW, DT, stable pairs) should follow.
Naively, we expect \cite{GV1, GV2, KKVSpinning} that within the class $\beta$,
the invariant $n_{r,\beta}$ counts curves of \emph{geometric} genus $r$ ``in
$X$". Here ``in" is to be interpreted loosely, including, as we discover below, maps
which are merely \emph{generically} embeddings. In
particular, a nonsingular curve $$C\subset X$$ of genus $g$ 
should contribute only
to $n_{g,\beta}$, while a reduced irreducible curve 
$C\subset X$ of arithmetic genus
$g$ and geometric genus $h<g$ should contribute to at most
$n_h,n_{h+1},\ldots,n_g$.

%We will be able to prove this (Theorem \ref{Kenny}) and the fact that,
%for a nodal curve, $n_k$ counts the number
%of finite maps from a curve of arithmetic genus $k$ to the class $\beta$
%that are generically injections (Proposition \ref{nodal}). That is a nodal
%curve $C\subset X$
%of arithmetic genus $g$ and geometric genus $g-r$ contributes to $n_k,\
%g\le k\le g-r$, exactly through each of its partial normalisations: normalising
%at any choice of $g-k$ nodes gives a finite, generically injective, map
%from a curve of arithmetic genus $k$ to $X$ with image $C$. Up to $\chi$-functions,
%reflecting the local deformation theory of $C$ and sheaves thereon, each
%such partial normalisation contributes $\pm1$ to $n_{k,C}$.

\subsection{Nonsingular curves}

The BPS contributions $n_{r,C}$ of a nonsingular curve $C$ are easy to compute.
First, we need to understand the local deformation theory of the pairs spaces
$P_n(X,\beta)$ about the locus $P_n(C)$ of pairs supported on $C$. The answer
turns out to be very simple with all of the $\chi^B$-functions $\chi\_P$
of these spaces being equal, up to sign, to the same constant $\chi\_\M(\O_C)$.
Here, as before, $\chi\_\M$ is the $\chi$-function of the moduli space
$$
\M_{1-g}(X,\beta)\ \ni\,\O_C
$$
of sheaves of Hilbert polynomial $k\mapsto k\int_\beta c_1(L)+(1-g)$.

\begin{lem} \label{om}
$\chi\_P|\_{P_n(C)}=(-1)^{n-1}\chi\_\M(\O_C)$ is constant.
\end{lem}

\begin{proof} We follow the proof of 
Theorem \ref{Kai}. As there, in a neighbourhood
of the locus
$$
\M_{1-g}(C)\subset\M_n(X,\beta)$$
of sheaves supported on $C$,
$\M_{1-g}(X,\beta)$ is isomorphic to all other such moduli spaces
$$
\M_{1-g+i}(X,\beta)\ \ni\,\O_C(D),
$$
where $D$ is a divisor on $C$ of degree $i$. We use the same argument
as before:
extend $D$ to a local divisor in $X$ (using
the nonsingularity of $C$) and map
$$
\spreaddiagramcolumns{1pc} \xymatrix{
\M_{1-g}(X,\beta) \rto^(.45){\otimes\O(D)} & \M_{1-g+i}(X,\beta),}
$$
taking $\O_C$ to $\O_C(D)$. Thus $$\chi\_\M(\O_C)=\chi\_\M(\O_C(D))$$ and the
$\chi^B$-function is \emph{identically constant} over the loci of sheaves
supported on $C$, since $D$ is arbitrary. By Theorem \ref{Kai},
the $\chi^B$-functions of all of the moduli spaces
$P_{1-g+i}(X,\beta)$ of pairs take the constant value $(-1)^{-g+i}\chi\_\M(\O_C)$
on restriction to $P_{1-g+i}(C)$.
\end{proof}

The space of pairs $P_{1-g+i}(C)$ supported on $C$ is the $i$th symmetric
product $S^iC$ of $C$ (a proof of a more general fact is given in Proposition
\ref{rhilb} below). By Lemma \ref{om}, the local stable pairs invariants
are
$$
P_{1-g+i,C} =(-1)^{-g+i}\chi\_\M(\O_C)e(S^iC).
$$

In this case the interpretation of the formula \eqref{localBPS} is
clear. For any manifold $M$, the generating function for the numbers $(-1)^ie(S^iM)$
is
\beq{eulergen}
(1+q)^{-e(M)}=
1-e(M)q+\frac{e(M)(e(M)+1)}2q^2-\ldots.
\eeq
Therefore $Z_{C}(q)$, the contribution of $C$ to $Z_{P,\beta}(q)$, is
\begin{multline*}
(-1)^g\chi\_\M(\O_C)\big(q^{1-g}-e(C)q^{2-g}+e(S^2C)q^{3-g}
-e(S^3C)q^{4-g}+\ldots\big) \\
=(-1)^g\chi\_\M(\O_C)q^{1-g}(1+q)^{2g-2}.
\end{multline*}
Since $q^{1-g}(1+q)^{2g-2}$
 is precisely the contribution of $n_{g,C}$ to \eqref{localBPS},
we have proved the following.

\begin{prop} A nonsingular curve $C\subset X$ of genus $g$ contributes
$$
n_{g,C}=(-1)^g\chi\_\M(\O_C)
$$
to $n_{g,\beta}$. And $n_{r,C}=0$ for $r\neq g$.
\end{prop}

\subsection{Singular curves: discussion}
For smooth curves, the geometry of the formulae \eqref{localBPS} is very
simple.  Remarkably,
the BPS formalism makes sense in the singular case also. 
To start,
we expand the formulae \eqref{localBPS}
out and read off the $n_{r,C}$ inductively:
\begin{eqnarray} \label{eqns}
P_{1-g,C} &=& n_{g,C}\,, \nonumber \\
\hspace{18mm} P_{2-g,C} &=& -e(\Sigma_g)\ n_{g,C}+n_{g-1,C}\,, \\ \nonumber
P_{3-g,C} &=& e(S^2\Sigma_g)\ n_{g,C}-e(\Sigma_{g-1})\ n_{g-1,C}+n_{g-2,C}\,,
\end{eqnarray}
and so on. Here and below we denote a smooth compact 2-manifold of genus
$g$ by $\Sigma_g$.

The formulae \eqref{eqns} tell us, inductively, what $C$ contributes
to each BPS number. The moduli space $P_{1-g}(C)$ consists of
the single point $(\O_C,1)$. By Theorem \ref{Kai},
 $P_{1-g,C}=(-1)^g\chi\_\M(\O_C)$, so by \eqref{eqns},
$$
n_{g,C}=(-1)^g\chi\_\M(\O_C).
$$
The contribution of the $n_{g,C}$ term to $P_{2-g,C}$ is then
\begin{equation} \label{jjyy}
-(-1)^ge(\Sigma_g)\chi\_\M(\O_C).
\end{equation}
If $C$ is nonsingular,
\eqref{jjyy} is precisely the contribution
of the space of pairs $P_{2-g}(C)\cong C$ supported on $C$, but for $C$
singular $P_{2-g,C}$ is the more complicated weighted Euler characteristic
$$
e(P_{2-g}(C),\chi\_P|\_{P_{2-g}(C)}).
$$
We \emph{define} $n_{g-1,C}$
to be the discrepancy between these two Euler characteristics:
\beq{d1}
n_{g-1,C}=P_{2-g,C}+(-1)^ge(\Sigma_g)\chi\_\M(\O_C).
\eeq

For example, consider
a curve $C\subset X$ with 1 node
for which the moduli space of sheaves $\M_{1-g}(X,\beta)$ is nonsingular
in a neighbourhood of $\M_{1-g}(C)$. Let $\pm$ denote
the sign $(-1)^{g+\dim\M_{1-g}(X,\beta)}$. Then $n_{g,C}=\pm1$, 
$P_{2-g}(C)\cong
C$, and so
$$
n_{g-1,C}=\pm(-e(C)+e(\Sigma_g))=\mp1.
$$
\smallskip
% \begin{eqnarray*}
% \tilde n_{g,\beta}(C) &=& 1, \quad\mathrm{and} \\
% \tilde n_{g-1,\beta}(C) &=& -e(C)+2-2g=-1.
% \end{eqnarray*}

We proceed inductively by viewing
 $C$ as contributing $n_{g,C}$ nonsingular curves of genus $g$ and
 $n_{g-1,C}$ nonsingular curves of genus $g-1$.
These genus $g$ and $g-1$ curves contribute
$$
e(S^2\Sigma_g)n_{g,C}-e(\Sigma_{g-1})n_{g-1,C}
$$
to $P_{3-g,\beta}$, as in \eqref{eqns}. The discrepancy
$$
n_{g-2,C}=P_{3-g,C}-e(S^2\Sigma_g)n_{g,C}+e(\Sigma_{g-1})n_{g-1,C}
$$
is what we define to be the number of genus $g-2$ curves supported
on $C\subset X$.

These formulae quickly become unmanageable, which is why we use the more
concise generating functions \eqref{localBPS}, to which they are equivalent.
A number of miraculous 
cancellations of Euler characteristics
and $\chi^B$-functions must occur for a singular curve
of geometric genus $h$ to have $n_{r,C}=0$ for $r<h$.
We will obtain these cancellations from an interplay between
Serre duality and Theorem \ref{Kai}.

\subsection{Singular curves: results}
Let $\Jac(C)$ denote the group of line bundles of degree 0 on $C$. There
is an action of $\Jac(C)$ on
 $\M_n(C)$ by tensoring. Let
$$\{ O_i \subset \M_n(C)\}_{i\in I}$$
be the set of $\Jac(C)$-orbits.
The index set $I$ need not be finite.
By convention we fix $O_0$ to be the orbit consisting of line bundles of
degree $n-1+g$ on $C$.

Fix a local effective divisor $D_C$ intersecting $C$ transversely at
a nonsingular point.
Tensoring with multiples of $\O(D_C)$ fixes local isomorphisms between
all of the $\M_n(X,\beta)$ in a neighbourhood of $\M_n(C)$. These isomorphisms
preserve the orbit types $O_i$. Therefore, we think of the $\M_n(X,\beta)$
and their subloci $O_i$ as (locally) independent of $n$.

The duality \eqref{duall} maps every orbit $O_i\subset\M_n(C)$ to another
orbit $O_{i'}\subset\M_{-n}(C)$ since
$$
\ext^2(F\otimes L,K_X)\cong L^{-1}\otimes\ext^2(F,K_X)
$$
for any (local) line bundle $L$. Hence, we obtain an involution
\beq{dualll}
i\leftrightarrow i'
\eeq
on the indexing set $I$.

% This has fixed points, orbits $O_i$ taken
%to themselves by \eqref{duall}, such as $O_0$ on any curve, or \emph{any}
%orbit $O_i$ on a nodal curve. 

\begin{lem} \label{Oi}
The restriction $\chi\_\M|\_{O_i}$ of the $\chi^B$-function of
$\M_n(X,\beta)$ is a constant $\chi\_i$ on each orbit $O_i$ (independent
of $n$). Moreover $\chi\_i=\chi\_{i'}$.
\end{lem}

\begin{proof} Every line bundle $L\in\Jac(C)$ 
can be trivialised over the finite
singular set of $C$, so is linearly equivalent to a (noneffective) divisor
supported on the nonsingular locus of $C$. The latter can 
be extended to a local divisor in $X$. Therefore $L$ is the restriction
of a line bundle $L$ defined on a neighbourhood of $C\subset X$. Then, the
map
$$
\xymatrix{\M_n(X,\beta) \rto^(.47){\otimes L} & \M_n(X,\beta),}
$$
defined only locally in a neighbourhood of $\M_n(C)$, is a local isomorphism
and so preserves $\chi^B$-functions. Thus,
$$
\chi\_\M(F)=\chi\_\M(F\otimes L)
$$
for any $F\in\M_n(C)\subset\M_n(X,\beta)$ and $L\in\Jac(C)$.

Since the construction commutes with the isomorphisms 
$$\M_n(C)\to\M_{n+1}(C)$$
given by tensoring with $\O(D_C)$, $\chi\_i$ is independent of $n$.
Since the duality \eqref{duall} preserves $\chi^B$-functions,
$\chi\_i= \chi\_{i'}$.
\end{proof}

Let $Z\subset C$ be the support of the singularities of $C$, 
and let $$C^0=C\take
Z$$ 
be the nonsingular locus.
Let $$\mu_Z=(2-2g)-e(C^0)$$
be  the Euler characteristic of the Milnor
fibre of $Z$: the sum over the components of $Z$ of one minus the Milnor
number of the component. The invariant $\mu_Z$ depends only on the analytic
germ of $Z\subset C$. Define 
$$
Q_n(C,Z)\subset P_{1-g+n}(C)
$$
to be the locus of pairs whose cokernel $Q$ (which has length $n$) is supported
entirely at $Z$.

Define an \emph{admissible} subset
\beq{admis}
J=\bigcup_nJ_n\subset\bigcup_n\M_n(C)
\eeq
to be a constructible
subset invariant under tensoring with line bundles on $C$ of any degree,
and invariant under the duality \eqref{duall}. Equivalently $J$ is a (possibly
infinite) union of orbit pairs $O_i\cup O_{i'}$ (if $i=i'$ is fixed
by the duality \eqref{dualll} then this is just ${O_i}$) and their translations
by line bundles of nonzero degree.

We set $Q_n(C,Z)_J$ to be the locus of points in $Q_n(C,Z)$ whose underlying
sheaf lies in $J$. 

\begin{prop} \label{Euler}
For any admissible $J$, the generating function of signed
topological Euler characteristics
$$
\sum_{n=0}^\infty(-1)^ne(Q_n(C,Z)_J)\ q^n
$$
can be written as
$$
\sum_{r=0}^gn_r(J)\ q^{g-r}(1+q)^{2r-2g-\mu_Z}
$$
for integers $n_r(J),\ r=0,\ldots,g$.
\end{prop}

\begin{proof}
We first use the same argument as in Proposition
\ref{ggtc}, applied to the
Euler characteristics of the fibres of the map $\phi_n$ from 
$$P_n(C)_J =\phi_n^{-1}(J)$$ to $J_n$. 
The fibres are empty for $n\le-g$. 

The fibre $\PP(H^0(F))$ over $F\in\M_n(C)$ has Euler characteristic
$h^0(F)$ while the fibre $\PP(H^1(F)^*)$ over the dual
$$
\ext^2_X(F,K_X)\in\M_{-n}(C)
$$
has Euler characteristic $h^1(F)$. If $F\in J_n$, then $\ext^2(F,K_X)\in
J_{-n}$. Since the difference between these two
Euler characteristics is $n$, we obtain
$$
e_n-e_{-n}=(-1)^{n+g-1}n\ e(J_1), \quad\text{and}\quad e_n=0,\ \ n\le-g,
$$
where we have defined
 $$e_n=(-1)^{n+g-1}e(P_n(C)_J)$$ and used 
 $J_n\cong
J_1$ 
for all $n$. Hence, we can write
\beq{beepees}
\sum_ne_nq^n=\sum_{r=0}^gn_r(J)\ q^{1-r}(1+q)^{2r-2}
\eeq
uniquely with
integers $n_r(J)$, $r=0,\ldots,g$, with $n_0(J)=(-1)^ge(J_1)$.

Restricted to the nonsingular locus $C^0$,
 every stable pair corresponds to a nonzero section
of a line bundle. The orbit type of the underlying sheaf is determined
by its restriction to the germ of $Z\subset C$. Therefore, 
stratifying $P_n(C)_J$
by the length of that part of the cokernel $Q=F/\O_C$ of the pair supported
on $Z$, we obtain the decomposition
\beq{decomp}
P_n(C)_J=\coprod_{k=0}^{n+g-1}Q_k(C,Z)_J\times S^{n+g-1-k}(C^0).
\eeq
Taking generating series of signed Euler characteristics gives, by \eqref{eulergen},
$$
\sum_n e_nq^n=
\sum_{k=0}^\infty(-1)^ke(Q_k(C,Z)_J)\ q^{k+1-g\,}(1+q)^{-e(C^0)}.
$$
Substituting \eqref{beepees} gives
$$
\sum_{r=0}^gn_r(J)q^{1-r}(1+q)^{2r-2}=q^{1-g}(1+q)^{2g-2+\mu_Z}
\sum_{k=0}^\infty(-1)^ke(Q_k(C,Z)_J)\,q^{k\!}.
$$
Rearranging implies the result.
\end{proof}

The advantage of the spaces $Q_k(C,Z)$ is that they depend only on the germ
of $Z\subset C$ (in the analytic topology), and not on $C$ or the genus $g$.

\begin{lem} \label{same}
Any curve with the same singularity germ as $Z\subset C$ has orbits $O_i,\,i\in
I$ and admissible subsets $J$  in 1-1 correspondence with
those of $C$. Moreover,
the spaces $Q_k(C,Z)$ and $Q_k(C,Z)_J$ and the duality
$i\leftrightarrow i'$ depend only on the germ of $Z\subset C$.
\end{lem}

\begin{proof}
Stable pairs supported on $C$ with given cokernel $Q$ are determined by the
extension
$$
0\to\O_C\to F\to Q\to0
$$
whose class lies in
$$
\Ext^1_C(Q,\O_C)=H^0(\ext^1_C(Q,\O_C)).
$$
The right side is just $\ext^1_C(Q,\O_C)$ thought of as a $\C$-module.
Thus the stable pair is determined by data entirely \emph{local to
the support of $Q$}.

Since the pairs $Q_k(C,Z)\subset P_{1-g+k}(C)$ have cokernel supported entirely
at $Z$, they are determined by the germ of $Z\subset C$.

The sheaf $F$ underlying a stable pair is a line bundle on $C^0$.
Hence, the 
orbit type is determined by the  restriction of $F$
to the germ of $Z\subset C$.
Similarly we claim the action of the duality \eqref{duall} on orbits $O_i$
is determined locally about $Z\subset C$ by
$$
F\mapsto\hom_C(F,\omega_C).
$$
Here, we consider $F$ as an $\O_C$-module and use relative Serre
duality for the embedding $\iota\colon C\into X$:
$$
R\hom_X(\iota_*F,K_X)\cong\iota_*R\hom_C(F,\omega_C)[-2].
$$
Taking $h^2$ of both sides and using the exactness of $\iota_*$ gives
$$
\ext^2_X(\iota_*F,K_X)\cong\iota_*\hom_C(F,\omega_C)
$$
as required.

Hence, the admissible subsets $J$ are also determined by the germ of $Z\subset
C$.
\end{proof}

We can now compare $C$ to a curve of minimal
genus with the same singularities. We let $g=g(C)$ denote the arithmetic
genus of $C$, and $g(\bar C)$ its geometric genus -- the genus of the
normalisation $$p\colon\bar C\to C.$$
 The \emph{$\delta$-invariant} of $Z$
is the difference:
$$
\delta(Z)=g(C)-g(\bar C)=\operatorname{length}\big((p_*\O_{\bar
C})\big/\O_C\big).
$$
The $\delta$-invariant
depends only on the germ of $Z\subset C$ (since $(p_*\O_{\bar C})/\O_C$
is supported at $Z$). It
is the minimal arithmetic genus of a
 curve $C'$ containing $Z$ as a singularity:
the arithmetic genus of a rational curve with singularity $Z$. Such
curves exist in Calabi-Yau 3-folds.

\begin{lem} \label{germ}
Given a germ of a singularity $Z\subset C$ of
a space curve $C$, there exists
a smooth projective 3-fold $X'$ containing a compact rational curve
$C'$ with singularities $Z'$ such that the germ of $Z'\subset C'$ is
isomorphic to the germ of $Z\subset C$.
We may also arrange that the canonical bundle $K_{X'}$ is trivial on all
curves close to $C'$.
\end{lem}

\begin{proof}
Normalise $C$ at $Z$,
$$
p\colon\bar C\to C,
$$
and take an analytic neighbourhood of $p^{-1}(Z)$ consisting of a finite
union of discs. Embed these discs
into $\PP^1$ and reapply $p$ to them to give a rational curve $C'$ with
singularities exactly $Z$.

The curve
$C'$ can be embedded in $\PP^N$ and then mapped to $\PP^3$ by a generic 
projection. Since the Zariski tangent spaces all have dimension $\le3$,
the map $C' \rightarrow \PP^3$
will still be an embedding.

Fix a nonsingular divisor $D\subset\PP^3$ of degree 2
 which intersects $C'$
transversally away from $Z$. Blow up $D\cap C'$ to give a 3-fold $X'$ with exceptional divisor $E$. The proper
transform $\bar C'\cong C'$ is a rational curve in $X'$.

All nearby curves are also the proper transforms $\bar C''$ of curves
$$C''\subset\PP^3$$ which pass through the points $D\cap C'$. 
For all such curves, $K_{\PP^3}|_{C''}$ is isomorphic to the divisor $-(2D\cap C')\subset C''$. Therefore,
$$
K_{X'}|_{\bar C''}\cong K_{\PP^3}|_{C''}+2E\cap C''=-2D\cap C'+2D\cap C'
$$
is trivial.
\end{proof}

Instead of studying $Q_k(C,Z)$ through $C$ of genus $g$, we can study
the same space $Q_k(C',Z')$ (by Lemma \ref{same}) through $C'$ of minimal
arithmetic genus $\delta(Z)$. By Proposition \ref{Euler},
we obtain
a stronger result:
\beq{Euler2}
\sum_{n=0}^\infty(-1)^ne(Q_n(C,Z)_J)q^n=
\sum_{r=0}^{\delta(Z)}n_r(J)\ q^{\delta(Z)-r}(1+q)^{2r-2\delta(Z)-\mu_Z}.
\eeq

\begin{thm} \label{Kenny}
The curve $C$ contributes only to the BPS numbers $n_r$ for $r$
between the geometric and arithmetic genera of $C$:
$$
n_{r,C}=0 \quad\text{for }r<g(\bar C)\text{ or }r>g(C).
$$
\end{thm}

\begin{proof}
For each $j\in\Z$ set
$$
J(j)=\bigcup_nJ(j)_n\subset\bigcup_n\M_n(C)
$$
to be the locus where the $\chi^B$-function $\chi\_\M$
restricted from $\M_n(X,\beta)$ equals $j$. By Lemma \ref{Oi},
$J(j)$ is invariant
under tensoring with line bundles and under the duality \eqref{duall}. Hence,
$J(j)$ is admissible \eqref{admis}: a (possibly infinite) union of orbit
pairs $O_i\cup
O_{i'}$ and their translates by line bundles of nonzero degree.

By \eqref{decomp} we have the stratification
$$
P_n(C)_{J(j)}=\coprod_{k=0}^{n+g-1}Q_k(C,Z)_{J(j)}\times S^{n+g-1-k}(C^0),
$$
where, by Lemma \ref{same}, the loci $Q_k(C,Z)_{J(j)}$ depend only on the
germ of $Z\subset C$. By
Theorem \ref{Kai} and \eqref{eulergen},  $P_n(C)_{J(j)}$ contributes
$$
\sum_{k=0}^\infty(-1)^{k-g}j\ e(Q_k(C,Z)_{J(j)})\ q^{k+1-g}(1+q)^{-e(C^0)}
$$
to $Z_{C}(q)$. By \eqref{Euler2}, the above equals
$$
(-1)^gj\sum_{r=0}^{\delta(Z)}n_r(J(j)) \ q^{\delta(Z)-r+1-g}
(1+q)^{2r-2\delta(Z)-\mu_Z-e(C^0)}.
$$
Setting $s=r+g(\bar C)$, we rewrite the contribution as
$$
(-1)^gj\sum_{s=g(\bar C)}^gn_{s-g(\bar C)}(J(j))q^{1-s}(1+q)^{2s-2}.
$$

Since $\chi\_\M$ takes only finitely many values, 
we add up over \emph{finitely
many} $J(j)$ to get
$$
Z_{C}(q)=\sum_{s=g(\bar C)}^gn_{s,C\,}q^{1-s}(1+q)^{2s-2},
$$
where $n_{s,C}=(-1)^g\sum_jjn_{s-g(\bar C)}(J(j))$.
We have obtained
 the required BPS form, with nonzero BPS numbers $n_{s,C}$ only for
$g(\bar C)\le s\le g$.
\end{proof}

\label{sr}
\subsection{Nodal curves}
We illustrate the results of Section \ref{sr} in the 
case of nodal curves, where the
formulae are rather simpler.

Consider first an irreducible elliptic curve $C$ with one node and
normalisation
$$
\PP^1\cong\bar C\Rt{p}C.
$$
The moduli space of sheaves $\M_1(C)$ is just a copy of $C$. There are two
orbits $O_i$ of $\Jac(C)$: $O_0$ is the smooth part $C^0$ of $C$ corresponding
to degree one line bundles on $C$, and $O_1$
is the nodal point $Z$ corresponding to the pure sheaf $p_*\O_{\bar C}$.

The corresponding constant $\chi\_\M$-functions are 
$$\chi\_0=\chi\_\M(\O_C), \ \ \chi\_1=\chi\_\M(p_*\O_{\bar C}).$$
 As we have seen,
$$
n_{1,C}=(-1)^g\chi\_0=-\chi\_0.
$$

Already for $n=1$ we have cohomology vanishing and $P_1(C)$ is a $\PP^0$-bundle
over $\M_1(C)\cong C$. Therefore,
$$
P_{1,C}=\chi\_0e(C^0)+\chi\_1e(Z)=\chi\_1.
$$
By \eqref{d1},
$$
n_{0,C}=P_{1,C}-0.\chi\_0=\chi\_1.
$$
By Theorem \ref{Kenny},
these
are the only
two BPS numbers: $n_{1,C}$ coming from $\O_C$ and $n_{0,C}$ coming from
$p_*\O_{\bar C}$.

To see the geometry involved in the formulae for the local invariants,
consider $P_2(C)$.
The moduli space $P_2(C)$ 
is a $\PP^1$-bundle over $\M_2(C)\cong C$.
Hence,
\beq{P2C}
P_{2,C}=-2P_{1,C}.
\eeq
Since $C$ is Gorenstein, we have the
identification
$$P_n(C)\cong\Hilb^n(C)$$
by Proposition \ref{RHilb} of Appendix B. 
Stratifying by both the number of points supported at the node $Z$
and by the type of the underlying sheaf (the subscript 0 denoting a line
bundle, 1 the
push-down of a line bundle from the normalisation), we expand the equality
\eqref{P2C} as
\begin{multline*}
-\chi\_0e(S^2C^0)-\chi\_1e(C^0)e(Z)-\chi\_0e(\Hilb^2(C,Z)_0)-
\chi\_1e(\Hilb^2(C,Z)_1) \\ =-2(\chi\_0e(C^0)+\chi\_1e(Z)).
\end{multline*}
We use $\Hilb^i(C,Z)$ to denote the subset of $\Hilb^iC$ supported at $Z$.
Since $C^0\cong\C^*$, $e(S^i(C^0))=0$ for all $i$, and $Z$ is a single
point. We obtain
$$
-\chi\_0e(\Hilb^2(C,Z)_0)-\chi\_1e(\Hilb^2(C,Z)_1)=-2\chi\_1.
$$
As in Section \ref{sr}, the equation can be considered as two equations holding
independently for $\chi\_0$ and $\chi\_1$ by applying Serre duality over
the appropriate orbit $O_i\subset\M_2(C)$ instead of the whole moduli space.
Then, we find
\beq{20}
e(\Hilb^2(C,Z)_0)=0 \quad\text{and}\quad e(\Hilb^2(C,Z)_1)=2.
\eeq
Similarly $P_n(C)$ is a $\PP^{n-1}$-bundle over $C$, so
$P_{n,C}=(-1)^{n-1}nP_{1,C}$. We obtain the following
generalization
of \eqref{20}:
\beq{i0}
e(\Hilb^n(C,Z)_0)=0 \quad\text{and}\quad e(\Hilb^n(C,Z)_1)=n
\quad\forall n\ge1.
\eeq

These consequences of Serre duality are what ensure there are no further
BPS numbers. In such a simple case, they can be verified directly.
Since $C$ has 2-dimensional Zariski tangent space at $Z$,
$$\Hilb^2(C,Z)\cong\PP^1,$$ with two distinguished
points $\{0,\infty\}$ corresponding to the directions of the two branches
of $C$ at $Z$. Then,
$$
\Hilb^2(C,Z)_0=\PP^1\take\{0,\infty\}\cong\C^*
\quad\text{and}\quad \Hilb^2(C,Z)_1=\{0,\infty\},
$$
implying \eqref{20}. Similarly $\Hilb^n(C,Z)_1$ is $n$ points, corresponding
to the $n$ different ways of distributing $(n-1)$ points between the two
points in $\bar C=\PP^1$ that lie over $Z$. Finally, 
since
$\Hilb^n(C,Z)_0$ admits a $\C^*$-action without fixed points, we
recover \eqref{i0}.

In general, it is much harder to obtain the required results
without Serre duality and the BPS generating
function formalism. Even when the singularity has a $\C^*$-action, localisation
is of limited use as it is hard to know on which orbit $O_i$ a fixed point
lies.
\medskip

Next, we consider curves with more nodes. Fix $C$ with arithmetic genus $g$,
geometric genus $g-r$, and $r$ nodes. 
Let $i=0,\ldots,2^{r-1}$ enumerate all
of the partial normalisations
$$
p_i\colon\bar C_i\to C
$$
of $C$ including the identity $p_0$. These also enumerate the orbits of
$\Jac(C)$ on the moduli spaces of sheaves $\M_n(C)$, where
$O_i$ corresponds
to push-downs $p_{i*}$ of line bundles from $C_i$.

Stable pairs can also be pushed down from a partial normalisation
$\bar C_i$. Given a stable pair on $\bar C_i$,
$$
\O_{\bar C_i}\stackrel{s}{\to} L,
$$
we can push down $L$ and $s$  and compose with the canonical section
of $p_{i*}\O_{\bar C_i}$,
$$
\O_C\to p_{i*}\O_{\bar C_i}\Rt{p_{i*}s}p_{i*}L,
$$
to give a stable pair on $C$. Let $g_i$ denote the arithmetic genus of
$\bar C_i$, and let $\bar C_i^0\subset \bar C_i$ be
the nonsingular locus. We see
$$
S^{n-1+g_i}(\bar C_i^0)\subset P_n(C)$$
by the push-down construction.

In fact, by \eqref{i0} applied to both $C$ and the $\bar C_i$, 
the only parts of the stratification (by $O_i$
orbit type) of the moduli space $P_n(C)$ of stable pairs that contribute
to $P_{n,C}$ are
\beq{decomp2}
S^{n-1+g}(C^0)\sqcup\coprod_iS^{n-1+g_i}(\bar C_i^0).
\eeq
The $i$th stratum contributes with 
constant $\chi^B$-function $(-1)^{n-1}\chi\_i$.
The Euler characteristic of the nonsingular locus $C^0$ equals
$e(\Sigma_g)$ and
$$
e(S^jC^0)=e(S^j\Sigma_g).
$$
Similarly
$$
e(S^j\bar C_i^0)=e(S^j\Sigma_{g_i}).
$$
Therefore, \eqref{decomp2} yields
$$
P_{n,C}=(-1)^{n-1}\left(\chi\_0e(S^{n-1+g}(\Sigma_g))+
\sum_i\chi\_ie(S^{n-1+g_i}(\Sigma_{g_i}))\right).
$$
Since the $i$th term is the same as the contribution of $\pm\chi\_i$ genus
$g_i$ nonsingular curves in $X$, we see 
$$
n_{h,C}=\sum_{i\colon g_i=h}(-1)^{g_i}\chi\_i.
$$

\begin{prop} \label{nodal}
$n_{h,C}=(-1)^r\sum_{i\colon g_i=h}\chi\_i.\hfill\square$
\end{prop}

Hence, the contributions of $C$ to $n_h$ come from all of the partial
normalisations of $C$. Proposition \ref{nodal}
 is closely related to an early interpretation of the
Gopakumar-Vafa invariants $n_h$ as the number of $h$-dimensional tori in
the moduli spaces $\M_n(C)$ of sheaves. The tori may be viewed to be
the orbits $O_i$,
the (push-downs of) Jacobians of the curves $\bar C_i$.

\appendix
\section{Reduced curve classes} \label{reduced}

The methods of Section \ref{one} do not extend to reduced classes $\beta$
since the sheaf $F$ underlying a stable pair $(F,s)$ need not be stable.
For instance, $F$ could be reducible with non-scalar automorphisms. Hence,
viewing stable pairs via the underlying sheaves decorated with a section
is no longer profitable.

However, the local approach of Section \ref{LocalSec}, viewing pairs via their
cokernels, works much the same
as in the irreducible case. In fact, since the singularities
$Z\subset C$ of any reduced
curve also appear, locally analytically, in an irreducible curve, the proof
of Theorem \ref{Kenny} goes through almost unchanged. We require  
the
analogue of Theorem \ref{Kai} for reduced reducible
curves and check the invariance of $\chi^B$ under
dualisation \eqref{duall}.
%Such results are most easily
%proved using the language of Artin stacks (\emph{defining} 
%the $\chi^B$-function
%of a local quotient $A/G$ to be $(-1)^{\dim G}$ times that of $A$), but since
%more elementary methods suffice we 
%sketch here the arguments.

\begin{thm} \label{Pconj}
Let $(F,s)$ be a stable pair on $X$, and let $L$ be a degree $\ell$
line bundle
on the support of $F$.
Then,
$$
\chi\_P(F,s)=(-1)^\ell \chi\_P(F\otimes L,t),\  \ \ \chi\_P(F,s)
=\chi\_P(\ext^2(F,K_X),t),
$$
whenever the sections $t$ on the right define stable pairs.
\end{thm}

\begin{proof}
The result
 can be proved for arbitrary stable pairs --- the details will appear in a
future paper \cite{ptnew}. 
For the application to reduced $\beta$, we only need the result
for stable pairs with \emph{reduced} support. So, we work locally about
a pair $(F,s)$ with support given by a reduced curve $C$ 
with $m$ irreducible components.

Just as in the proof of Theorem \ref{Kai}, there exists a smooth
local divisor $D_C$ which intersects each irreducible component of $C$ transversally
in a single smooth point of $C$, disjoint from the zeros of $s$. Hence, $D_C$
also intersects all nearby curves $C'$ in $m$ smooth points.

Let $k\gg0$. Tensoring sheaves $F$ with $\O(kD_C)$ and multiplying their
sections $s$ by the canonical section $s_{D_C}^k$ of $\O(kD_C)$ yields a
local embedding,
\beq{miracle2} \xymatrix{
P_n(X,\beta)\supset\hspace{-1cm} & V_n\ \ \ar@{^{(}->}[r]^(.4){\phi} &
\ \ V_{n+mk} & \hspace{-1cm} \subset P_{n+mk}(X,\beta)\,.}
\eeq
Here, $V_r$ is a sufficiently small analytic neighbourhood of the locus $P_r(C)\subset
P_r(X,\beta)$ of stable pairs supported on $C$.
 
By making $V_{n+mk}$ smaller if necessary, $\phi$ admits a left inverse
$$ \xymatrix{
V_n & V_{n+mk} \lto_(.55){\psi}}
$$
given by forgetting the $mk$ points close to $D_C$. 
%More precisely,
%for any pair close to the image of $\phi$, we remove the part of its cokernel
%supported close to $D_C$ (and this part has length $mk$).
 More precisely,
in a neighbourhood of $D_C$ we replace any pair with support $C'$ by the
trivial pair $(\O_{C'},1)$. Since these are smooth points of $C'$,
 the map
$\psi$ is locally a smooth fibre bundle with fibre the $(mk)$-th symmetric
product of an open set of a smooth curve. 

By  relation \eqref{smooth}, 
$\phi$ and $\psi$ simply multiply $\chi^B$-functions
by $(-1)^{mk}$. Therefore, replacing  $(F,s)$ by
 the image $(F(kD_C),s_i\cdot
s^k{D_C})$ under $\phi$, we may assume that
 \emph{$F$ has no higher cohomology}. \medskip

We first explain how to prove Theorem \ref{Pconj} 
in case $L$ is trivial.
Let $Q_n(X,\beta)$ denote the moduli space of triples $(F,s_1,s_2)$, where

\vspace{6pt}

\begin{itemize}
\item $F$ is a pure sheaf with 1-dimensional support and Hilbert polynomial
$\chi(F(k))=\int_\beta c_1(\O(k))+n$, 
\item $s_1,s_2\in H^0(F)$ such that $\xymatrix{\O_X\oplus\O_X \rto^(.66){(s_1,s_2)}
& F}$ has finite cokernel,
\end{itemize}

\vspace{6pt}

\noindent
modulo the obvious equivalence induced by automorphisms of $F$. While 
the corresponding moduli functor is easily seen to be 
separated and
proper (and these triples have no automorphisms), the moduli space does
not seem to have been constructed by GIT
in the literature, 
though some very similar moduli problems have been
treated in \cite{HLFramed, LPPairs2}.
Recently Schmitt \cite{Schmitt} has constructed just such a moduli
space of \emph{torsion free} sheaves, and his methods certainly 
extend to the
\emph{pure} sheaves above. Further
 details will appear in \cite{ptnew}.

Let $U\subset Q_n(X,\beta)$ denote the open set on which both sections $s_i$
have finite cokernel. There are two projections from $U$ to $P_n(X,\beta)$,
mapping $(F,s_1,s_2)$ to $(F,s_1)$ and $(F,s_2)$.
Both maps are smooth with fibre
an open set in $H^0(F)$ of dimension $n$ since 
$H^1(F)=0$. By Proposition
1.5(i) of \cite{BehrendDT},
$$
\chi\_P(F,s_1)=(-1)^n\chi\_Q(F,s_1,s_2)=\chi\_P(F,s_2),
$$
as required.
\medskip

Consider now 
$F\otimes L$ for some degree 0 line bundle $L$ defined in a neighbourhood
of the support $C$ of the $F$. 
Since $L$ can be trivialised over any finite
number of points of $C$, we can write
$$
L\cong\O(D_1-D_2)
$$
for effective locally defined Cartier divisors $D_i$ intersecting $C$ transversally
in points disjoint from the singularities of $C$ and the cokernel of
the sections $s_i$. Let $d$ denote the degree of $D_1$ (and $D_2$)
on $C$.

By the same working as in \eqref{miracle2}, the $\chi^B$-functions of
\beq{1}
(F,s) \quad\text{and}\quad (F(D_1+kD_C),s\cdot s\_{D_1}\cdot s_{D_C}^k)
\eeq
differ only by $(-1)^{d+mk}$. Similarly for
\beq{2}
(F\otimes L,t) \quad\text{and}\quad (F\otimes L(D_2+kD_C),t\cdot
s\_{D_2}\cdot s_{D_C}^k).
\eeq
By construction,
$$
F(D_1+kD_C)=F\otimes L(D_2+kD_C).
$$
Then, 
the $\chi^B$-functions of \eqref{1} and \eqref{2}
agree by the previous proven case.
\medskip

If $F\otimes L$ where the degree $L$ is nonzero, we
apply \eqref{miracle2} separately to both $F$ and $F\otimes L$
and reduce to
the degree 0 case.
\medskip

Finally, consider $\ext^2(F,K_X)$. 
Here, we form a moduli
space of triples $(F,s_1,s_2)$, where
\vspace{6pt}

\begin{itemize}
\item $F$ is a pure sheaf with 1-dimensional support and Hilbert polynomial
$\chi(F(k))=\int_\beta c_1(\O(k))+n$,
\item $s_1\in H^0(F(r))$ has finite cokernel,
\item $s_2\in H^0(     \ext^2(F,K_X) (r))$ has finite cokernel,
\end{itemize}

\vspace{6pt}
\noindent and $r$ is a fixed integer.

The
corresponding problem
 is separated
(but not proper) and also has a quasi-projective
moduli space which has never been explicitly written down in the literature.
The moduli space
 maps to $P_{r+n}(X,\beta)$ by forgetting $s_2$
and to $P_{r-n}(X,\beta)$ by forgetting $s_1$. These maps are
smooth fibrations of relative dimensions $r-n$ and $r+n$ 
respectively for $r\gg0$. Hence,
the $\chi^B$-function of the triple $(F,s_1,s_2)$ equals both
$$
(-1)^{r-n}\chi\_P(F(r),s_1) \quad\text{and}
\quad (-1)^{r+n}\chi\_P(   \ext^2(F,K_X)(r),s_2).
$$
The required result is then obtained from previous cases.
\end{proof}

The proof of Theorem \ref{Kenny} now applies verbatim to give an expression
\beq{ZCq}
Z_{C}(q)=\sum_{s=g(\bar C)}^gm_{s,C\,}q^{1-s}(1+q)^{2s-2},
\eeq
as before, proving Theorems 1 and 2 for reduced curve classes.

However, there is difference now.
For a reducible curve, $$g(\bar
C)=1-e(\bar C)/2$$ might be negative. 
Here, the Gopakumar-Vafa conjecture
takes the form of a sum over decompositions of $C$ into
reducible components{\footnote{The $C_i$ have maximal arithmetic genus $g_i$, but need not be irreducible.
We do \emph{not} sum over different orderings of the $C_i$.}}:
$$
Z_C(q)=\sum_{\sum_iC_i=C}\ 
\prod_i\ \sum_{r_i=0}^{g_i}\Big(n_{r_i,C_i}\,q^{1-r_i}(1+q)^{2r_i-2}\Big).
$$
To properly define the BPS counts of $C$,
 we would need to show that all
of the negative genus contributions of \eqref{ZCq} can be written inductively
as products of BPS counts of lower degree curves. 
We currently have only partial results in this direction.

\section{Stable pairs on Gorenstein curves} \label{B}
\subsection{Gorenstein curves}
The scheme theoretic support $C$ of a stable pair $(F,s)$
is always Cohen-Macaulay. Equivalently,
the
dualising complex of $C$
is a \emph{sheaf} $\omega_C$. If $\omega_C$
is a \emph{line bundle}, then the curve is \emph{Gorenstein}.
Plane curves are basic examples. If $C$ can be embedded
in a nonsingular surface $S$, then 
$$\omega_C\cong
\omega_S(C)|_C$$
where $\omega_S$ is the canonical line bundle of $S$.
Hence, $C$ is Gorenstein.

Let $(F,s)\in P_{1-g+n}(C)$ be a stable pair supported on 
a Gorenstein curve $C$. Dualising
$$
\O_C\Rt{s}F
$$
on $C$ yields the map
\beq{pairhilb}
F^*\To\O_C.
\eeq
We will show \eqref{pairhilb} is the
ideal sheaf of a length $n$ subscheme $Z\subset C$.
Conversely, we  associate a stable pair to each subscheme
$[Z] \in \Hilb^n(C)$ and
establish a bijection
$$P_{1-g+n}(C) \leftrightarrow \Hilb^n(C).$$
In Appendix \ref{sru}, the corresponding isomorphism
of {\em schemes} is established for curves $C$ in surfaces.

\begin{lem} Let $C$ be a Gorenstein curve, and let $F$ be
a sheaf on $C$ which is generically locally free. Then, $F$ is
pure if and only if
$$\ext^i_C(F,\O_C)=0
\quad \forall  i > 0\ . $$  \label{pureext}
\end{lem}

\begin{proof}
Suppose $F$ is pure.
Since $F$ is generically locally free, 
the sheaf $\ext^i_C(F,\O_C)$ is supported
at a finite number of points. 
Hence, the vanishing 
\begin{equation}\label{ggww8}
H^0(\ext^i_C(F,\O_C)\otimes\mathcal L)=0 \quad \forall i > 0 
\end{equation}
for any line bundle $\mathcal L$ on $C$ implies the Lemma.

Using the Gorenstein condition on $C$, 
let $\mathcal L=\omega_C\otimes L $. 
For $L\gg0$, by the local-to-global spectral sequence and vanishing,
\eqref{ggww8} is isomorphic to
$$
\Ext^i_C(F\otimes L^*,\omega_C)\cong H^{1-i}(F\otimes L^*)^*,
$$
using Serre duality. But $L^*\ll0$ and $F$ is pure, so 
$H^0(F\otimes L^{-1})=0$. Thus $H^{1-i}(F\otimes L^{-1})=0$ for all $i>0$.

Conversely, suppose $F$ is generically locally free on $C$
 with vanishing $\ext^1(F,\O_C)$. Consider the sequence
\beq{KFE}
0\to K\to F\to E\to0
\eeq
where $K$ is the largest subsheaf of $F$ with $0$-dimensional support. 
The quotient $E$ is 
pure and generically locally free. 
We have proven above the vanishing
$$\ext^i(E,\O_C)=0\quad \forall i>0\ . $$ 
By \eqref{KFE}, we conclude
$$
\ext^1(K,\O_C)=\ext^1(F,\O_C),
$$
which vanishes by assumption. Hence,
$$
0=H^0(\ext^1(K,\O_C))=\Ext^1(K\otimes\omega_C,\omega_C)=H^0(K\otimes\omega_C)^*
$$
by Serre duality. Since $K$ has 0-dimensional support, 
$K$ must vanish.
Therefore, $F$ is pure.
\end{proof}

\begin{prop} \label{rhilb}
A stable pair supported on a 
Gorenstein curve $C$ is equivalent to a 0-dimensional
subscheme of $C$. Under the equivalence, the pair
\begin{equation} \label{veffg}
0 \rightarrow \O_C \stackrel{s}{\rightarrow} F \rightarrow Q 
\rightarrow 0 
\end{equation}
is associated to the subscheme 
$$\O_C \cong \ext^0(\O_C,\O_C)\rightarrow  \ext^1(Q,\O_C) \rightarrow 0\ .$$
\end{prop}

\begin{proof}
Given a stable pair $(F,s)\in P_{1-g+n}(C)$, 
we apply $\hom_C(\ \cdot\ ,\O_C)$ to the sequence \eqref{veffg}.
By the purity of $\O_C$ and
 Lemma \ref{pureext}, we obtain the exact sequence
$$
0\to F^*\to\O_C\to\ext^1(Q,\O_C)\to0.
$$
Hence, $F^*$ is an ideal sheaf and 
$\ext^1(Q,\O_C)$ is isomorphic to the structure sheaf
of a subscheme of $C$.
The higher terms in the sequence
yield the vanishing of
 $\ext^{\ge2}(Q,\O_C)$, so $\ext^1(Q,\O_C)$ has length $n$.
We have defined a map of sets from
$P_{1-g+n}(C)$ to $\Hilb^n(C)$.

Given a subscheme $[Z]\in \Hilb^n(C)$, the dual of the exact sequence
$$
0\to\I_Z\to\O_C\to\O_Z\to0
$$
starts as
\beq{ubc}
0\to\O_C\to\I_Z^*\to\ext^1(\O_Z,\O_C)\to0.
\eeq
By Lemma \ref{pureext},
the higher terms in the sequence
yield the vanishing of
 $\ext^{\ge2}(O_Z,\O_C)$, so $\ext^1(\O_Z,\O_C)$ has length $n$
and
$\I_Z^*$ has holomorphic Euler characteristic $1-g+n$.

Since $\ext^{\geq 1}(\I_Z,\O_C)$ vanishes by Lemma \ref{pureext}, 
$\I_Z^*=R\hom(\I_Z,\O_C)$. By
applying $R\hom(\ \cdot\ ,\O_C)$ again, we obtain 
$$
R\hom(\I_Z^*,\O_C)=\I_Z.
$$
Therefore, $\ext^{\ge1}(\I_Z^*,\O_C)$ vanishes and
$\I_Z^*$ is pure by  Lemma \ref{pureext}. We conclude
$$  \O_C \rightarrow \I_Z^*$$
determines a stable pair in $P_{1-g+n}(C)$.
We have defined a map of sets from $\Hilb^n(C)$ to
$P_{1-g+n}(C)$.

As the two constructions are easily seen to be inverse
to each other, a bijection is established. 
\end{proof}

Simple examples where the quotient $Q$ of \eqref{veffg}
fails to be a structure sheaf can be found on nodal
curves. However, 
$\ext^1(Q,\O_C)$ is always a structure sheaf. For
the equivalence of Proposition \ref{rhilb}, the duals are
necessary.

The equivalence of Proposition \ref{rhilb} is continuous.
For Gorenstein curves,
the moduli of pairs $P_{1-g+n}(C)$ is 
homeomorphic as
a topological space 
to the Hilbert scheme $\Hilb^n(C)$. 
If $C$ is also \emph{reduced}, we can apply the result of
Proposition \ref{Euler}, and the method of Theorem \ref{Kenny}, replacing
the $\chi^B$-function by $(-1)^\mathrm{dim}$.
The result is that the generating function of Euler characteristics
$$
\sum_{i=0}^\infty e(\Hilb^i(C))q^i
$$
can be written in the form
$$
\sum_{r=g(\bar C)}^g n_{r,C\,}q^{g-r}(1-q)^{2r-2}
$$
for integers $n_{r,C}$, where $r$ runs from the geometric genus $g(\bar C)$
to the arithmetic genus $g=g(C)$ of $C$.

\subsection{Surfaces} \label{sru}
Let $S$ be a nonsingular projective surface, and let
 $\beta\in H_2(S,\Z)$. Let $\M$ denote
the moduli space of pure dimension 1 subschemes of $S$ in class 
$\beta$. Since the holomorphic Euler characteristic of such
subschemes is determined by adjunction, $\M$ is a Hilbert scheme
of curves.
Let
$$\CC\subset S\times\M$$
 denote the universal curve.

By Proposition \ref{rhilb}, $P_{1-g+n}(S,\beta)$ is in bijective
correspondence with
the relative Hilbert scheme $\Hilb^n(\CC/\M)$, where
$$
2g-2=\int_S K_S\cdot \beta +\beta^2 \ \ .
$$
We now show the correspondence is an isomorphism of schemes.

\begin{prop} \label{RHilb}
The moduli space of pairs $P_{1-g+n}(S,\beta)$ is isomorphic to the
relative Hilbert scheme $\Hilb^n(\CC/\M)$.
\end{prop}

\begin{proof}
Let $P$ denote 
the space of pairs $P_{1-g+n}(S,\beta)$, and let $(\FF,s)$ denote the
universal stable pair on $S\times P$. 
The kernel of the canonical map
$$
\O_{S\times P}\Rt{s}\FF
$$
is the ideal sheaf of a relative curve $\mathcal D$ flat over $P$. 
By the universal property of $\M$,
 we obtain
a map
 $$f\colon P\to\M$$
such that
 $\mathcal D=(1_S\times f)^*\CC$ is the pull-back of the universal curve.

We have the canonical exact sequence
$$
0\to\O_{\mathcal D}\Rt{s}\FF\to\Q\to0,
$$
with the universal quotient $\Q$ also flat over $P$. 
Applying $\hom_{\mathcal D}(\ \cdot\ ,\O_{\mathcal D})$ yields
\beq{relseq}
0\to\FF^*\to\O_{\mathcal D}\to\ext^1(\Q,\O_{\mathcal D})
\to\ext^1(\FF,\O_{\mathcal D}).
\eeq
The last term vanishes by flatness, base change, and 
Lemma \ref{pureext}. 
Similarly $\ext^{\ge2}(\FF,\O_{\mathcal D})=0$, 
which implies $\ext^{\ge2}(\Q,\O_{\mathcal D})=0$. 
Therefore $\ext^1(\Q,\O_{\mathcal D})$ is flat over $P$ with relative length $n$.

By the universal property of $\Hilb^n(\CC/\M)$, 
the quotient sequence \eqref{relseq} defines a map from $P$ 
to $\Hilb^n(\CC/\M)$ through which $f$ factors. 
The inverse map is similarly a relative version of the map in the proof of Proposition \ref{rhilb}.
\end{proof}

\section{$K3$ surfaces}
\subsection{Nonsingularity}
Interesting irreducible classes on Calabi-Yau 3-folds can be found
on $K3$ fibrations. Let
$$\pi:X \rightarrow \bigtriangleup$$
be a fibration of a Calabi-Yau 3-fold over a nonsingular curve satisfying
the following properties:
\begin{enumerate}
\item[(i)] the fibres of $\pi$ are $K3$ surfaces,
\item[(ii)] $S=\pi^{-1}(0)$ carries a
            irreducible (1,1)-class $\beta\in H_2(S,\mathbb{Z})$,
\item[(iii)] $\pi$ is transverse to the Noether-Lefschetz
             locus associated to $\beta$.
\end{enumerate}
For the above geometry, we do not require $X$ to be compact.

Condition (iii) is equivalent to requiring the
Kodaira-Spencer class
$$\kappa\in H^1(S,TS) \cong H^{1,1}(S)$$
to evaluate to something nonzero on $\beta$.
By (iii), a curve on $S$ in class $\beta$ does not deform,
even to first order, away from $S$.
Therefore,
\begin{equation} \label{jjw}
P_n(S,\beta)\rightarrow P_n(X,\iota_*\beta)
\end{equation}
is a bijective correspondence with
a component of the moduli
space of pairs, where
$$\iota:S \rightarrow X.$$
We will show \eqref{jjw} is a local isomorphism of schemes
in Lemma \ref{zzzz} below.

The moduli space $P_n(S,\beta)$ is empty if $n< 1-g$, 
where
$$2g-2 =\int_S \beta^2 $$
as before.
           
\begin{prop} \label{nonsing}
For $\beta$ irreducible, \label{nmw}
the moduli space $P_n(S,\beta)$ is
nonsingular of dimension $n+2g-1$. 
\end{prop}

\begin{proof}
We briefly recall the deformation theory of pairs on
surfaces \cite{HePairs,LPPairs1,
LPPairs2}. The deformations of the pair $\O_S\to F$ are governed by
$\Hom(I\udot,F)$
and the obstruction theory by $\Ext^1(I\udot,F)$. The maps
\beq{trmap}
\Ext^1(I\udot,F)\to\Ext^2(F,F)\Rt{\tr}H^2(\O_S)
\eeq
take the obstructions to deforming $\O_S\to F$ first to the obstructions
to deforming the sheaf $F$ and then \cite{MuSymplectic} to the obstruction to deforming
the determinant $\O_S(C)$
of the sheaf $F$, where $C$ is the support of $F$.
Since the Picard scheme is smooth, the latter obstruction vanishes.
Hence, the obstructions to deforming
$F$ can be taken to lie in the trace-free group $\Ext^2(F,F)_0$, and
the obstructions to deforming the pair $\O_S\to F$ lie in the
kernel
of the map 
$$\Ext^1(I\udot,F)\to H^2(\O_S)$$ 
obtained from \eqref{trmap}. If the kernel is 0,
then by the analysis in \cite{HePairs,LPPairs1, LPPairs2}, the moduli
space
is nonsingular.

The deformation and obstruction spaces sit inside
the exact sequence
\begin{multline} \label{sseq}
0\to\Hom(F,F)\to H^0(F)\to\Hom(I\udot,F)\to\Ext^1(F,F)\to \\
H^1(F)\to\Ext^1(I\udot,F)\to\Ext^2(F,F)\to0,
\end{multline}
induced by the triangle $I\udot\to\O_S\to F$.
We claim that the first arrow on the second line  is zero,
or equivalently, that
$$
\Ext^1(F,F)\to H^1(F)
$$
is onto. It is enough to show the composition
\begin{equation}\label{nnw}
H^1(\O_C)\Rt{\id}H^1(\hom(F,F))\subseteq\Ext^1(F,F)\to H^1(F)
\end{equation}
is onto. But \eqref{nnw} is multiplication by the section $s$, sitting in the
exact
sequence
$$
H^1(\O_C)\rt{s}H^1(F)\to H^1(Q).
$$
Since the support of $Q$ is 0-dimensional, $H^1(Q)$ vanishes.

The support $C$ of $F$ is reduced and irreducible since
 $C$ is Cohen-Macaulay and $\beta$ is irreducible. 
Since the
rank of $F$ on $C$ is 1, $F$ is
simple. Therefore the identity map
$$
H^0(\O_S)\Rt{\id}\Hom(F,F)
$$
is an isomorphism. By Serre duality and $K_S\cong\O_S$, the trace
map
$$
\Ext^2(F,F)\Rt{\tr}H^2(\O_S)
$$
is also an isomorphism.

So \eqref{sseq} has become the sequences
$$
0\to H^0(F)\big/\langle s\rangle\to\Hom(I\udot,F)\to\Ext^1(F,F)\to
H^1(F)\to0
$$
and
\beq{ext20}
\Ext^1(I\udot,F)\Rt{\sim}\Ext^2(F,F)\cong H^2(\O_S)=\C.
\eeq
Therefore, the obstruction space 
$$\ker\left(\Ext^1(I\udot,F)\to
H^2(\O_S)\right)$$
vanishes, and the moduli space is nonsingular.

The dimension can be easily computed. The space of curves of
class $\beta$ has dimension $g$. For a nonsingular curve $C$,
the dimension of the space of pairs $P_n(C)$ is
$n+g-1$. Hence, the dimension of $P_n(S,\beta)$ is
$n+2g-1$.
\end{proof}

Proposition \ref{nonsing} was first proven by Kawai-Yoshioka in \cite{ky},
where the space of stable pairs on a surface is interpreted as a moduli
space of D2-D0 branes. We include the above argument for completeness here.
With hindsight, one can see many of the techniques we use 
in 3-dimensions
in the 2-dimensional analysis of \cite{ky}.

\begin{lem}  \label{zzzz}
For $\beta$ irreducible,
$P_n(S,\beta)\subset P_n(X,\iota_*\beta)$
is isomorphic to a component of the moduli space of pairs.
\end{lem}

\begin{proof}
Let $I\udot_X=\{\O_X\to\iota_*F\}$ and $I\udot_S=\{\O_S\to F\}$. There
is a canonical exact triangle
$$
F(-S)\to L\iota^*(I\udot_X)\to I\udot_S
$$
on $S$.
Applying $\Hom_S(\ \cdot\ ,F)$ gives
\begin{multline*}
0\to\Hom_S(I\udot_S,F)\to\Hom_X(I\udot_X,\iota_*F)\to\Hom_S(F,F(S))
\\ \to\Ext^1_S(I\udot_S,F)\to\Ext^1_X(I\udot_X,\iota_*F)\to\ldots\ \ .
\end{multline*}

From the simplicity of $F$, we see
$$\Hom_S(F,F(S)) \cong H^0(\O_S(S)).$$
Using the identification \eqref{ext20},
we can rewrite the above long exact sequence as
$$
0\to\Hom_S(I\udot_S,F)\to\Hom_X(I\udot_X,\iota_*F)\to H^0(\O_S(S))\to
H^2(\O_S)\to\ldots
$$
The map from $H^0(\O_S(S))\cong T_0\bigtriangleup$ to
$H^2(\O_S)\cong\C$ takes a normal direction to the fibre 
$S\subset X$ to the obstruction to deforming the determinant of
$I\udot$ sideways.
The determinant is $\O_S(-C)$, where $C$ is the support of $F$, and the
obstruction to deforming sideways is $\langle\kappa,\beta\rangle\ne0$. 
Therefore
$$
\Hom_S(I\udot_S,F)\cong\Hom_X(I\udot_X,\iota_*F).
$$

All the deformations of the pair $(\iota_*F,\iota_*s)$ are the push-forwards
of deformations of the pair $(F,s)$. By Proposition \ref{nmw},
the moduli space of the latter is nonsingular. Hence, the
moduli space of the former is also nonsingular and
the push-forward map is a local isomorphism.
\end{proof}

\subsection{BPS states}

Let $P_n(S,h)$ denote the moduli space for an irreducible
class $\beta$ satisfying
\beq{betah}
2h-2= \int_S \beta^2.
\eeq
Let $\Omega_{P}$ be the cotangent bundle of
the moduli space $P_n(S,h)$.
The self-dual 
obstruction theory on $P_n(S,h)$ induced from
the inclusion \eqref{jjw} has obstruction bundle
$\Omega_{P}$. Hence, the contribution of $P_n(S,h)$
to the stable pairs invariants of $X$ is
$$Z^S_h(y) = \sum_{n} (-1)^{n+2h-1} e(P_n(S,h))\  y^n.$$

Fortunately, the topological Euler characteristics of 
$P_n(S,h)$ have been calculated by Kawai-Yoshioka.
By Theorem 5.80 of \cite{ky},
\begin{multline*}
\sum_{h=0}^\infty \sum_{n=1-h}^\infty  
e(P_n(S,h))\  y^n q^h = \\
\left(\sqrt{y}-\frac{1}{\sqrt{y}}\right)^{-2}\
\prod_{n=1}^\infty \frac{1}{(1-q^n)^{20} (1-yq^n)^2(1-y^{-1}q^n)^2} \ .
\end{multline*}
For our pairs invariants, we require the signed Euler characteristics,
$$
\sum_{h=0}^\infty Z_h^S(y) \ q^h = 
\sum_{h=0}^\infty \sum_{n=1-h}^\infty (-1)^{n+2h-1} 
e(P_n(S,h))\  y^n q^h. 
$$
Therefore,
$\sum_{h=0}^\infty Z_h^S(y) \ q^h$ is
$$
-\left(\sqrt{-y}-\frac{1}{\sqrt{-y}}\right)^{-2}\
\prod_{n=1}^\infty \frac{1}{(1-q^n)^{20} (1+yq^n)^2(1+y^{-1}q^n)^2} \ .
$$

Let $r_{g,h}$ be the BPS invariant in genus $g$ and class $\beta$
satisfying \eqref{betah}.
By the definition of the
BPS invariants for the theory of stable
pairs \eqref{PT1BPS},
$$\sum_{h=0}^\infty Z_h^S(y) \ q^h = \sum_{g=0}^\infty \sum_{h=0}^\infty
(-1)^{g-1} r_{g,h}
\left(\sqrt{-y}-\frac{1}{\sqrt{-y}}\right)^{2g-2} q^h.$$
Putting all the formulae together yields the following result.
\begin{prop}\label{kkww}
We have
\begin{multline*}
\sum_{g=0}^\infty \sum_{h=0}^\infty
(-1)^{g} r_{g,h}
\left(\sqrt{z}-\frac{1}{\sqrt{z}}\right)^{2g} q^h = \\
\prod_{n=1}^\infty \frac{1}{(1-q^n)^{20} (1-zq^n)^2(1-z^{-1}q^n)^2}\ . 
\end{multline*}
\end{prop}
Proposition \ref{kkww} is exactly the Katz-Klemm-Vafa 
\cite{KKVSpinning} prediction
for BPS state counts for irreducible classes on a $K3$ surface.
Our proof is really just
an interpretation of the calculation of Kawai-Yoshioka in
the theory of stable pairs.

If $\beta\in H_2(S,\mathbb{Z})$ is primitive, then the geometry
can be deformed to make $\beta$ irreducible, and Proposition \ref{kkww}
still applies. If $\beta$ is {\em not} primitive, the Katz-Klemm-Vafa
formula prediction is not yet proven.
On the Gromov-Witten side, the Katz-Klemm-Vafa formula
is open even in the primitive case. See \cite{dm} for a discussion.

\subsection{Yau-Zaslow}
We can specialize the Katz-Klemm-Vafa formula of Proposition \ref{kkww}
to genus 0 by letting $z\rightarrow 1$,
\begin{equation}
 \sum_{h=0}^\infty
 r_{0,h} q^h = 
\prod_{n=1}^\infty {(1-q^n)^{-24}}\ . 
\end{equation}

Consider the linear system of curves on $S$ of irreducible
class $\beta$ satisfying
$$2h-2=\int_S \beta^2.$$
By Theorem \ref{Kenny}, only curves of geometric genus 0 contribute
to the BPS count $r_{0,h}$. By the main result of \cite{xc},
the only genus 0 curves on a generic $K3$ (with algebraic
class $\beta$) are {\em nodal}.
By Proposition \ref{nodal}, a nodal rational curve $C$  contributes
$$(-1)^0 \sum_{i\colon g_i=0} \chi\_i = (-1)^0(-1)^{n-1} \chi^B_{P_n(C)}
=(-1)^{n-1}(-1)^{n+2g-1} =1$$ 
to $r_{0,h}$.

We conclude the strong enumerative form of the Yau-Zaslow formula:
$r_{0,h}$ exactly counts rational curves on a generic $K3$
surface with $\beta$ algebraic.

\addtocontents{toc}{\SkipTocEntry}

\vspace{+14 pt}

\begin{minipage}[position]{65mm}
\noindent
Departement Mathematik\\
ETH Z\"urich\\
{\tt rahul@math.ethz.ch}
\end{minipage}
\begin{minipage}[position]{8cm}
\noindent Department of Mathematics\\
Imperial College\\
{\tt rpwt@imperial.ac.uk}
\end{minipage}

\end{document}